\documentclass{amsart}

\usepackage{amsmath}
\usepackage{amssymb}
\usepackage{graphicx}
\usepackage{color}
\usepackage[all]{xy}
\usepackage{amsthm}
\usepackage[colorlinks,linkcolor=blue,citecolor=blue,pdfstartview=FitH]{hyperref}

\theoremstyle{plain}
\newtheorem{thm}{Theorem}
\newtheorem{prop}{Proposition}[section]
\newtheorem{lemma}[prop]{Lemma}
\newtheorem{cor}[prop]{Corollary}

\newtheorem{remark}{Remark}
\newtheorem{remarks}{Remarks}

\theoremstyle{remark}

\theoremstyle{definition}
\newtheorem{example}{Example}
\newtheorem{dfn}[prop]{Definition}
\newtheorem{dfns}[prop]{Definitions}

\newcommand{\wt}[1]{\widetilde{#1}}
\newcommand{\co}{\colon\thinspace}

\newcommand{\bound}{\partial}
\newcommand{\ph}{\varphi}

\newcommand{\Q}{\mathbb{Q}}
\newcommand{\Z}{\mathbb{Z}}
\newcommand{\C}{\mathbb{C}}
\newcommand{\tr}{\mathrm{Tr}}
\newcommand{\rk}{\mathrm{Rk}}
\newcommand{\pie}{\pi_1}
\renewcommand{\O}{\mathcal{O}}
\newcommand{\m}{\mathfrak{m}}

\newcommand{\calJ}{\mathcal{J}}
\newcommand{\calI}{\mathcal{I}}

\newcommand{\calG}{\mathcal{G}}

\newcommand{\calB}{\mathcal{B}}
\newcommand{\calS}{\mathcal{S}}
\newcommand{\calM}{\mathcal{M}}
\newcommand{\sfX}{{\sf X}}
\newcommand{\sfY}{{\sf Y}}
\newcommand{\sfR}{{\sf R}}

\newcommand{\sfXA}{{\sf X}_A}

\renewcommand{\sl}{\text{SL}_2}
\newcommand{\psl}{\text{PSL}_2}

\begin{document}

\title{Closed surfaces and character varieties}

\author{Eric Chesebro}
\address{Department of Mathematical Sciences, University of Montana} 
\email{Eric.Chesebro@mso.umt.edu} 

\begin{abstract}  The powerful character variety techniques of Culler and Shalen can be used to find essential surfaces in knot manifolds.  We show that module structures on the coordinate ring of the character variety can be used to identify detected boundary slopes as well as when closed surfaces are detected.   This approach also yields new number theoretic invariants for the character varieties of knot manifolds.
\end{abstract}

\maketitle

\section{Introduction}

Suppose that $N$ is a compact, irreducible 3-manifold with torus boundary and that $\sfX$ is an irreducible algebraic component of the $\sl \C$-character variety for $N$.  It is known that the dimension of $\sfX$ is at least one \cite{CCGLS}.  It is also known that $\sfX$ contains a great deal of topological information about $N$.   

Let $\wt{\sfX} \to \sfX$ be a birational map from a smooth projective curve.  This map is defined for all but finitely many points.  These points are called ideal points.  For $\gamma \in \pie N$ and $\chi \in \sfX$, define $I_\gamma(\chi)=\chi(\gamma)$.  This determines a rational function $I_\gamma \co \wt{\sfX} \to \C$.  This paper concerns the following landmark result of Culler and Shalen.

\begin{thm}[Culler-Shalen] \label{thm:CS} If $\hat{x}$ is an ideal point, there is an associated non-empty essential surface $\Sigma$ in $N$.  
\begin{enumerate}
\item $\Sigma$ may be chosen to have empty boundary if and only if $I_\alpha$ is regular at $\hat{x}$ for every peripheral element $\alpha \in \pie N$.
\item Otherwise, there exists a unique slope with the property that if $\alpha \in \pie N$ represents this slope then $I_\alpha$ is regular at $\hat{x}$.  In this case, every component of $\bound \Sigma$  represents the slope corresponding to $\alpha$.
\end{enumerate}
\end{thm} 
When an essential surface arises from $\sfX$ and this theorem, we say that the surface is detected by $\sfX$.

Since 1983, case (2) of Theorem \ref{thm:CS} has been carefully studied.  Many famous papers have provided both applications and insight for this case.  Notably, the Culler-Shalen norm of \cite{CGLS} and the A-polynomial of \cite{CCGLS} can tell us exactly which slopes are detected.  These tools are remarkably effective and have been used in proofs of several famous theorems including the Smith Conjecture \cite{Sh2} and the Cyclic Surgery Theorem \cite{CGLS}.  In contrast, there is very little in the literature concerning case (1) of Theorem \ref{thm:CS}.

 It has also been difficult to interpret global properties of $\sfX$ in terms of the topology of $N$.  This paper introduces new functions which relate detected boundary slopes and essential closed surfaces in the manifold to module structures of the coordinate ring $\C[\sfX]$.   In the spirit of \cite{CL2}, these functions also contain number theoretic information about $\sfX$. 

The trace ring $T(\sfX)$ is the subring of the coordinate ring $\C[\sfX]$ which is generated by $\{ I_\gamma \, | \, \gamma \in \pie N \}$.  Define $T_\Q(\sfX)$ to be the smallest $\Q$-algebra in $\C[\sfX]$ which contains $T(\sfX)$.  Suppose that $\alpha \in \pie N$ is primitive and peripheral.  Then $\C[\sfX]$, $T_\Q(\sfX)$, and $T(\sfX)$ have the structure of a $\C[I_\alpha]$-module, a $\Q[I_\alpha]$-module, and a $\Z[I_\alpha]$-module, respectively.  Let $\rk_\sfX^\C(\alpha)$, $\rk_\sfX^\Q(\alpha)$, and $\rk_\sfX^\Z(\alpha)$ denote the ranks of these modules.  Note that the functions $I_\alpha$ are well-defined on slopes.  Hence, these rank functions may be viewed as functions on $\mathcal{S}$, the set of slopes for $N$.  It is clear that \begin{equation} \label{eq:ineq} \rk_\sfX^\C \ \leq \ \rk_\sfX^\Q \ \leq \ \rk_\sfX^\Z.\end{equation}
The following theorem is the main theorem of this paper.

\begin{thm} \label{thm:intro1}
Let $\mathbb{X} \in \{ \C, \Q, \Z\}$.  Then
\begin{enumerate}
\item The function $\rk_\sfX^\mathbb{X} \co \mathcal{S} \to \Z^+ \cup \{ \infty \}$ is constant with value $\infty$ if and only if $\sfX$ detects a closed essential surface.
\item Otherwise, $\rk_\sfX^\mathbb{X} (\alpha)=\infty$ if and only if $\sfX$ detects the slope $\alpha$.
\end{enumerate}
\end{thm}

We begin by proving the theorem for $\mathbb{X}=\C$.  A straightforward application of this case yields a proposition of independent interest. 

\begin{prop}
Let $\sfX^\circ$ be the union of the irreducible components $\sfX'$ of the $\sl \C$ character variety such that \begin{enumerate} \item $\sfX'$ detects a closed essential surface if and only if $\sfX$ does, and
\item for each slope $\alpha$, $\sfX'$ detects a surface with boundary slope $\alpha$ if and only if $\sfX$ does.
\end{enumerate}
Then $\sfX^\circ$ is defined over $\Q$.
\end{prop} 

This proposition, together with Theorem \ref{thm:CS}, also allows us to improve the result to include $\mathbb{X}=\Q$.  As an application of this result, we prove another proposition.

\begin{prop}
If $\sfX$ contains the character of a non-integral representation which detects a closed essential surface, then $\sfX$ detects a closed essential surface.
\end{prop}

Together with data compiled by Goodman, Heard, and Hodgson, this proposition gives many examples of hyperbolic knot manifolds for which the hyperbolic components of their character varieties detect closed essential surfaces.  We also use this proposition in our proof that the main theorem holds in the case $\mathbb{X}=\Z$.

In addition to our main result, we begin an investigation of the basic properties of the functions $\rk_\sfX^\mathbb{X}$ and include computations for several examples.  It is easy to prove the following two propositions.

\begin{prop} \label{prop:introX_A}
Suppose that $N$ is a knot manifold and $H_1(N;\Z) \cong \Z$, generated by $\alpha$.  Let $\sfXA \subseteq \sfX(N)$ be the curve consisting of all abelian representations of $\pie N$. Then  \[  \rk_{\sfX_A}^\C(\alpha) \ = \ \rk_{\sfX_A}^\Z(\alpha) \ =\ 1.\]  
\end{prop}

\begin{prop}
If $N$ is a knot manifold and $\sfX$ is a norm curve component of $\sfX(N)$ then $\rk_\sfX^\C(\alpha) \geq 2$ for every slope $\alpha$.
\end{prop}

As a corollary to Theorem \ref{thm:intro1}, we have the following.
\begin{cor}
Suppose $N$ is the exterior of a two-bridge knot and $\langle \mu, \beta \, | \, \omega \beta=\mu \omega \rangle$ is the standard presentation for $\pie N$.  If $\sfX \subset \sfX(N)$ is an irreducible algebraic component defined over $\Q$ then $\sfX$ is defined by an irreducible polynomial of the form \[ I_{\mu \beta}^n - \sum_{j=0}^{n-1} p_j(I_\mu) I_{\mu \beta}^j.\]  The set $\{ I_{\mu \beta}^j \}_0^{n-1}$ is a free basis for $\C[X]$ as a $\C[I_\mu]$-module, $T_\Q(\sfX)$ as a $\Q[I_\mu]$-module, and $T(\sfX)$ as a $\Z[I_\mu]$-module.
\end{cor}

We point out that, with the appropriate definitions, these results also hold for the $\psl \C$-character variety.  In what follows, we use $\sfX$ to indicate that we are working in the $\sl \C$ setting and $\sfY$ in the $\psl \C$ setting.  When $N \subset S^3$ is a knot exterior, we reserve $\lambda,\mu \in \pie N$ to be a longitude, meridian pair.  We write $\sfX_A$ and $\sfY_A$ for the algebraic sets of abelian characters.   When $N$ is hyperbolic, we use the notation $\sfX_0$ and $\sfY_0$ to indicate algebraic components of the character varieties which contain a discrete faithful character.  

Below, we list the results of our calculations in five different examples.  The last example is of particular interest.  It shows that the inequalities (\ref{eq:ineq}) need not be equalities and that, although $T(\sfX)$ must always be torsion free as a $\Z[I_\sigma]$-module, it need not be a free module.

%%%%%%%%%%%%%%%%%%%%%%%%%%%%%%%%%%%
\smallskip
\begin{enumerate}
\item Suppose $N$ is the exterior of a trefoil knot.  Then $\sfX(N)=\sfX_A \cup \sfX_0$ where $\sfX_0$ is irreducible. We have \[ \rk_{\sfX_0}^\Z(\mu) \ =  \ 1.\]  This shows that the converse to Proposition \ref{prop:introX_A} does not hold.

%%%%%%%%%%%%%%%%%%%%%%%%%%%%%%%%%%%
\smallskip
\item  Suppose $N$ is the exterior of the figure-eight knot.  We have 
\begin{gather*} 
 \rk_{\sfX_0}^\Z(\mu) \ = \ \rk_{\sfY_0}^\Z(\mu) \ = \ 2 \\
 \rk_{\sfY_0}^\C(\lambda) \ = \ \rk_{\sfY_0}^\Z(\lambda) \ = \ 4 \\
 \rk_{\sfY_0}^\C(\mu^2 \lambda) \ = \ \rk_{\sfY_0}^\Z(\mu^2 \lambda) \ = \ 5. 
\end{gather*}

%%%%%%%%%%%%%%%%%%%%%%%%%%%%%%%%%%%

\item Suppose $N=M_{003}$.  Then $\pie N=\langle  \gamma , \eta  \, | \,  \gamma  \eta  \gamma ^{-2} \eta  \gamma  \eta ^3 \rangle$.  Define $\mu=( \eta ^2 \gamma  \eta  \gamma )^{-1}$ and $\lambda=( \gamma  \eta  \gamma )^{-1} \eta  \gamma  \eta $.  Then $\mu$ and $\lambda$ are primitive, peripheral, and generate a peripheral subgroup of $N$.  We have
\[ \rk_{\sfX_0}^\C(\mu) \ = \ \rk_{\sfX_0}^\Z(\mu) \ =\  4.\]  Also, $I_\mu \in T(\sfY_0)$ and \[\rk_{\sfY_0}^\C(\mu) \ = \ \rk_{\sfY_0}^\Z(\mu) \ =\  2.\] Note that, since $\sfY_0$ is a norm curve, these ranks achieve their minimum possible value.  In contrast to examples (1) and (2), the $\psl \C$-rank is strictly smaller than the $\sl \C$-rank (at $\mu$).  As with examples (1) and (2),  all $\Z[I_\mu]$-modules are free.

%%%%%%%%%%%%%%%%%%%%%%%%%%%%%%%%%%%

\begin{figure}
\setlength{\unitlength}{.1in}
\begin{picture}(23,15)
\put(0,0) {\includegraphics[width= 2.3in]{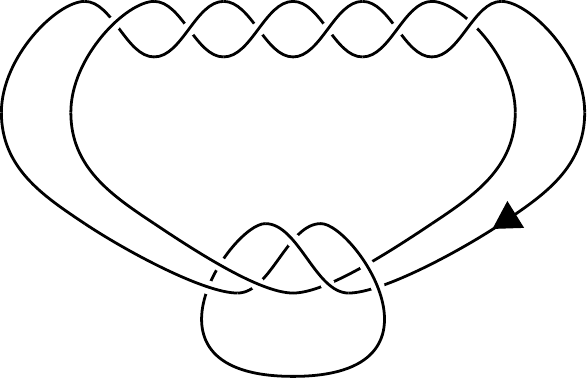}}
\put(1,5.5){$\sigma$}
\put(16, 1){$\mu$}
\put(4, 8.5){$\eta$}
\end{picture}
\caption{The knot $8_{20}$ labelled with Wirtinger generators}  \label{fig:8_20}
\end{figure}
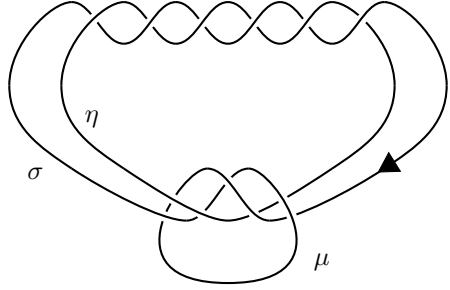

\smallskip
\item Suppose $N$ is the exterior of the knot $8_{20}$.  $N$ is hyperbolic, so we consider $\sfX_0$.  The diagram in Figure \ref{fig:8_20} gives a Wirtinger presentation for $\pie N$ which reduces to $\left\langle \mu, \gamma \, | \,  \mu  \gamma ( \mu  \gamma  \mu )^{-1} \gamma ^3( \mu  \gamma  \mu )^{-1} \gamma  \mu  \gamma ^{-2}\right\rangle$  where $\gamma=\eta \sigma$.  

The functions $I_\mu$, $I_\gamma$, and $I_{\gamma \mu}$ give an embedding of $\sfX_0$ into $\C^3$.  $I_{\gamma}$ and $I_{\gamma \mu}$ both satisfy irreducible integral dependencies in $\Z[I_\mu][x]$ of degree 5.  It follows that $\calB=\{ I_\gamma^i I_{\gamma \mu}^j \, | \, 0 \leq i,j \leq 4 \}$ is a generating set for $\C[\sfX_0]$ as a $\C[I_\mu]$-module, $T_\Q(\sfX_0)$ as a $\Q[I_\mu]$-module, and $T(\sfX_0)$ as a $\Z[I_\mu]$-module.  Notice that $| \calB |=25$.

A Groebner basis argument shows that relations amongst the elements of $\calB$ are plentiful and, in fact, $\{ 1,I_\gamma, I_\gamma^2, I_\gamma^3, I_\gamma^2 I_{\gamma \mu}, I_{\gamma \mu} \}$ is a free basis for each of the above modules.  Hence, \[ \rk_{\sfX_0}^\C(\mu) \ = \ \rk_{\sfX_0}^\Z(\mu) \ = \ 6\] and $T(\sfX_0)$ is a free $\Z[I_\mu]$-module.

%%%%%%%%%%%%%%%%%%%%%%%%%%%%%%%%%%%
\smallskip
\item Suppose $N$ is the once punctured hyperbolic torus bundle from Section 5 of \cite{D2}.  Then $\pie N = \langle \alpha,\beta,\tau \, | \, \tau\alpha\tau^{-1}=(\beta\alpha\beta)^{-1}, \tau\beta\tau^{-1}=\beta\alpha(\beta\alpha\beta)^{-3} \rangle$.  The elements $\tau$ and $\lambda=[\alpha,\beta]$ form a basis for the peripheral subgroup of $N$ and the functions 
\begin{align*} 
  t&=I_\tau &  u&=I_{\alpha\tau} & v&=I_{\beta\tau}  &  w&=I_{\alpha\beta\tau} \\
  x&=I_\alpha & y&=I_\beta & z&=I_{\alpha\beta} 
\end{align*}
 give an embedding of $\sfX(N)$ into $\C^7$.  For $\epsilon \in \{ 0,1 \}$, we have an irreducible algebraic component $\sfX_\epsilon \subset \sfX(N)$ which contains a discrete faithful character.  
 \begin{enumerate}
\item $\{ 1,z,z^2,z^3,y,zy, z^2y, z^3y\} $ is a free basis for $\C[\sfX_\epsilon]$ as a $\C[I_\lambda]$-module.  Hence, $\rk_{\sfX_\epsilon}^\C(\lambda)=8$.
\item $\{ 1,z,z^2,z^3,y,zy, z^2y, z^3y, t, zt, z^2t, z^3t, yt, zyt, z^2yt, z^3yt\}$ is a free basis for $T(\sfX_\epsilon)$ as a $\Q[I_\lambda]$-module.  Hence, $\rk_{\sfX_\epsilon}^\Q(\lambda)=16$.
\item $\{ 1, y, y^2, y^3, t, u, v, w, x, z, yt, yu, yx, y^2x, xt, xu, vy \} $ generates $T(\sfX_\epsilon)$ as a $\Z[I_\lambda]$-module.   However, $T(\sfX_\epsilon)$ is torsion free but \underline{not} free as a $\Z[I_\lambda]$-module.  In particular, $\rk_{\sfX_\epsilon}^\Z(\lambda)=17$.
\end{enumerate}
The projection $(t,y,z)\co \sfX_0 \cup \sfX_1 \to \C^3$ is an isomorphism onto its image and the relation \[\rk_{\sfX_\epsilon}^\Q(\lambda) < \rk_{\sfX_\epsilon}^\Z(\lambda)\] reflects the fact that the inverse of this isomorphism is not defined over $\Z$.  The equality \[\rk_{\sfX_\epsilon}^\Q(\lambda)=2\cdot \rk_{\sfX_\epsilon}^\C(\lambda)\] reflects the fact that $\sfX_\epsilon$ is not defined over $\Q$.  
\end{enumerate}

\subsection{Acknowledgements}  We thank Steve Boyer, Brendan Hasset, and Kelly McKinnie for helpful conversations and suggestions.

%%%%%%%%%%%%%%%%%%%%%%%%%%%%%%%%%%%

%%%%%%%%%%%%%%%%%%%%%%%%%%%%%%%%%%%

\section{Algebraic Geometry} \label{sec:geometry}

Throughout this section we take $k$ to be an algebraically closed field.  All varieties are defined over $k$.  When $\sfX$ is an irreducible variety we write $k[\sfX]$ and $k(\sfX)$ to denote the ring of regular functions on $\sfX$ and the function field for $\sfX$, respectively.   If $p \in \sfX$, define $\O_{\sfX,p}$ to be the ring of germs of functions which are regular on neighborhoods of $p$.   We have a surjective homomorphism $k[\sfX] \to k$ given by $f \mapsto f(p)$.  Let $\m_p$ be the maximal ideal which is the kernel of this map.  By Theorem 3.2 part (c) of \cite{H}, $\O_{\sfX,p}$ is isomorphic to the localization $k[\sfX]_{\m_p}$.  We use this isomorphism to identify these two rings and think of $\O_{\sfX,p}$ as a subring of $k(\sfX)$.

\begin{dfns} Suppose $A \subseteq B$ are commutative rings and $1_B \in A$. 
\[\] \vspace{-.5in}
\begin{enumerate}
\item  An element $b \in B$ is {\it integral over $A$} if there is a number $n \in \Z^+$ and $\{ \alpha_i \}_0^{n-1} \subset A$ such that \[ b^n + \alpha_{n-1} b^{n-1} + \dots + \alpha_0 \ = \ 0.\]
This equation is called an {\it integral dependence relation of $b$ over $A$}.
\item {\it $B$ is integral over $A$} if every element of $B$ is integral over $A$.
\smallskip
\item The {\it integral closure $\overline{A}$ of $A$ in $B$} is the set of all elements of $B$ which are integral over $A$.  
\smallskip
\item If $\overline{A}=A$ then $A$ is {\it integrally closed in $B$}.
\end{enumerate}
\end{dfns}

By Corollary 5.3 of \cite{AM}, the integral closure of $A$ in $B$ is a ring.  The following characterization of integral closures is useful.

\begin{thm}[Corollary 5.22 of \cite{AM}] \label{thm:intersect}
Suppose $A$ is a subring of a field $K$.  The integral closure of $A$ in $K$ is the intersection of all valuation rings of $K$ that contain $A$.
\end{thm}

\begin{dfns}  Suppose that $\sfX$ is an affine variety.
\begin{enumerate}
\item $\sfX$ is a {\it normal variety} if $k[\sfX]$ is integrally closed in $k(\sfX)$.
\smallskip
\item A {\it normalization of $\sfX$} is an irreducible normal variety $\sfX^\xi$ and a regular birational map $\xi \co \sfX^\xi \to \sfX$ where $k[\sfX^\xi]$ is integral over $\xi^\ast \left( k[\sfX] \right)$.
\end{enumerate}
\end{dfns}

If $\sfX$ is an affine variety and $p \in \sfX$ then the quotient $\m_p/\m_p^2$ has the structure of a $k$-vector space.  

\begin{dfns}  Suppose that $\sfX$ is an affine variety.
\begin{enumerate}
\item A point $p \in \sfX$ is {\it non-singular} if $\dim\left( \m_p/\m^2_p \right)$ equals the dimension of $\sfX$.
\item $\sfX$ is {\it non-singular} if every point in $\sfX$ is non-singular.
\end{enumerate}
\end{dfns}

Theorem 1 in Chapter II, Section 5.1 of \cite{S} states that non-singular affine varieties are normal.

It is clear that if $\m_p$ is principal then $\dim \left( \m_p/\m^2_p \right) =1$.  Theorem 2 in Chapter II, Section 5.1 of \cite{S} implies that if $\sfX$ is a normal affine algebraic curve and $p \in \sfX$ then $\m_p$ is principal.  Thus, we have the following well-known theorem.

\begin{thm} \label{thm:smoothmodel}
Suppose $\sfX$ is an irreducible affine algebraic curve.  $\sfX$ is non-singular if and only if $\sfX$ is normal.
\end{thm}

The following theorem is a direct consequence of Theorems 4 and 5 in Chapter II, Section 5.2 in \cite{S} (and the first part of the proof of Theorem 4).

\begin{thm} \label{thm:normalization}
If $\sfX$ is an irreducible affine algebraic curve then $\sfX$ has a normalization $\xi \co \sfX^\xi \to \sfX$.  Moreover, the normalization is unique, affine, and its coordinate ring is the integral closure of $\xi^\ast\left(k[\sfX]\right)$.
\end{thm}

The projective coordinates on projective space $\mathbb{P}^n$ give distinguished open affine subsets $\{ U_i\}_0^n$ which cover $\mathbb{P}^n$.  So if $\sfX\subset \mathbb{P}^n$ is an irreducible projective variety then we have the distinguished affine open subsets $U_i \cap \sfX$ which cover $\sfX$.  We say that $\sfX$ is non-singular if $U_i \cap \sfX$ is non-singular for every $i = 0, \ldots , n$.

\begin{dfn} Suppose that $\sfX$ is an irreducible affine curve and let $\xi \co \sfX^\xi \to \sfX$ the normalization of $\sfX$.  As in \cite[Chapter 6]{H}, there is a smooth projective curve $\wt{\sfX}$ (unique up to isomorphism) so that $\sfX^\xi$ is isomorphic to an open set in $\wt{\sfX}$.  The projective variety $\wt{\sfX}$ is called the {\it smooth projective model for $\sfX$}.  
\end{dfn}

Henceforth, we identify $\sfX^\xi$ with its image in $\wt{\sfX}$ and we let $\iota \co \wt{\sfX} \dashrightarrow \sfX^\xi$ be the rational map which is the identity on $\sfX^\xi$.

\begin{dfn}
The {\it ideal points} of $\sfX$ are the points in the set $\calI(\sfX) = \wt{\sfX} - \sfX^\xi$.  It follows from the uniqueness of normalizations and smooth projective models that $\calI(\sfX)$ is well-defined.
\end{dfn}

We have dominant birational maps 
\[ \xymatrix{
\wt{\sfX} \ar@{.>}[r]^{\iota \ } &  \sfX^\xi \ar[r]^\xi & \sfX
} \]
The induced maps 
\[ \xymatrix{
k\big(\wt{\sfX}\big) &  k\left(\sfX^\xi\right) \ar[l]_{\ \, \iota^\ast} & k(\sfX)  \ar[l]_{\ \ \ \xi^\ast}
} \]
are isomorphims.   Moreover, $\xi^\ast \big( k[\sfX] \big) \subseteq k[\sfX^\xi]$.

Now suppose that $Y$ is an affine variey and $\ph \co \sfX \to  Y$ is a regular map with $\overline{\ph(\sfX)}=Y$.  Then $\ph^\ast \co k(Y) \to k(\sfX)$ is a field monomorphism and $\ph^\ast (k[Y]) \subseteq k[\sfX]$.

\begin{dfn}
A {\it hole} in $\varphi \co \sfX \to Y$ is a point $\hat{p} \in \calI(\sfX)$ such that 
$\left( \varphi \xi \iota\right)^\ast \left( k[Y] \right) \subseteq \O_{\wt{\sfX}, \hat{p}}$.
\end{dfn}

\begin{remark} If $\hat{p}$ is a hole in $\varphi \co \sfX \to Y$ then, since $\hat{p} \notin \sfX^\xi$, there is an element $f \in k[\sfX]$ with $(\xi \iota)^\ast (f) \notin \O_{\wt{\sfX}, \hat{p}}$ .    Also, if $\ph$ is a birational map, then $\ph$ is surjective if and only if $\ph$ has no hole.  \end{remark}

\begin{lemma} \label{lem:easy}
If $\varphi \co \sfX \to Y$ has a hole then $k[\sfX]$ is not integral over $\varphi^\ast \left( k[Y] \right)$.
\end{lemma} 
\proof Take a distinguished open affine set $U_i \subset \mathbb{P}^n$ with $\hat{p} \in U_i$.  Set $\sfX_{\hat{p}} = U_i \cap \wt{\sfX}$.   By Theorem \ref{thm:smoothmodel}, $\sfX_{\hat{p}}$ is a non-singular affine curve so $k[\sfX_{\hat{p}}]$ is integrally closed.  Recall that $\O_{\sfX_{\hat{p}},\hat{p}}$ is isomorphic to the localization $k[\sfX_{\hat{p}}]_{\m_{\hat{p}}}$.  Proposition 5.13 of \cite{AM} gives that this localization is integrally closed.  The inclusion induced isomorphism $k(\wt{\sfX}) \to k(\sfX_{\hat{p}})$ restricts to an isomorphism $\O_{\wt{\sfX}, \hat{p}} \to \O_{\sfX_{\hat{p}},\hat{p}}$.  Hence $\O_{\wt{\sfX},\hat{p}}$ is integrally closed.

The point $\hat{p}$ is a hole in $\varphi \co \sfX \to Y$ so $(\varphi \xi \iota)^\ast (k[Y]) \subseteq \O_{\wt{\sfX},\hat{p}}$.  Since $\O_{\wt{\sfX},\hat{p}}$ is integrally closed, every element of $k(\wt{\sfX})$ which is integral over $(\varphi \xi \iota)^\ast (k[Y])$ is an element of $\O_{\wt{\sfX},\hat{p}}$.  We have $f \in k[\sfX]$ with $(\xi \iota)^\ast (f) \notin \O_{\wt{\sfX},\hat{p}}$, hence $(\xi \iota)^\ast(k[\sfX])$ is not integral over $(\varphi \xi \iota)^\ast (k[Y])$ and so $k[\sfX]$ is not integral over $\varphi^\ast (k[Y])$.
\endproof

\begin{thm} \label{thm:geom}
Suppose $\sfX$ is an irreducible affine algebraic curve and $\varphi \co \sfX \to Y$ is a regular map with $\overline{\varphi(\sfX)}=Y$.  Then $\varphi$ has no hole if and only if $k[\sfX]$ is integral over $\ph^\ast \left( k[Y] \right)$.
\end{thm}

\proof 
By Lemma \ref{lem:easy}, we need only show that if $\ph$ has no hole then $k[\sfX]$ is integral over $\ph^\ast(k[Y])$.

Assume, to the contrary, that $k[\sfX]$ is not integral over $\ph^\ast(k[Y])$. Since $\xi^\ast$ is injective, this implies that  $k[\sfX^\xi]$ is not integral over $(\ph \xi)^\ast (k[Y])$.  Taking integral closures in $k(\sfX^\xi)$ we have
\[ k[\sfX^\xi] \ =\ \overline{k[\sfX^\xi]} \ \supsetneq \ \overline{(\ph \xi)^\ast (k[Y])}.\]
Take $f \in k[\sfX^\xi] - \overline{(\ph \xi)^\ast (k[Y])}$.  By Theorem \ref{thm:intersect}, there is a valuation ring $R \subset k(\sfX^\xi)$ with $f \notin R$ and $k \subset \overline{(\ph \xi)^\ast (k[Y])} \subseteq R$.

\smallskip
\noindent \underline{Claim:}   
$\iota^\ast (R) = \O_{\wt{\sfX},\hat{p}}$ for some $\hat{p} \in \wt{\sfX}$.

\smallskip
Using the claim, $\iota^\ast (f) \notin \O_{\wt{\sfX},\hat{p}}$ so $\hat{p} \in \calI(\sfX)$.   Moreover, $(\ph \xi)^\ast k[Y] \subseteq R$ so $(\ph \xi \iota)^\ast (k[Y]) \subseteq \iota^\ast (R) = \O_{\wt{\sfX}, \hat{p}}$.  That is, $\hat{p}$ is a hole in $\ph$, a contradiction.

The claim follows from Corollary 6.6 of \cite{H} since it gives an open set $\sfX_R$ in $\wt{\sfX}$ and a point $\hat{p} \in \sfX_R$ such that $\iota^\ast(R) = \O_{\sfX_R, \hat{p}} = \O_{\wt{\sfX}, \hat{p}}$. 
\endproof

It is well-known that, in this setting, $k[\sfX]$ is integral over $f^\ast (k[Y])$ if and only if $k[\sfX]$ is a finitely generated $f^\ast(k[Y])$-module (see for example Proposition 5.1 and Corollary 5.2 of \cite{AM}).   Hence, we have the following immediate corollary.

\begin{cor} \label{cor:geom} Suppose $\sfX$ is an irreducible affine algebraic curve and $\ph \co \sfX \to Y$ is a regular map with $\overline{\ph(\sfX)}=Y$.  Then $\ph$ has no hole if and only if $k[\sfX]$ is a finitely generated $\ph^\ast(k[Y])$-module.
\end{cor}

\begin{remarks} Suppose that $k[\sfX]$ is a finitely generated $\ph^\ast(k[Y])$-module.
%\[\] \vspace{-.4in}
\begin{enumerate}
\item Since $\sfX$ is irreducible, $k[\sfX]$ has no zero divisors.  Hence $k[\sfX]$ is a torsion free $\ph^\ast(k[Y])$-module.
\item We can be more concrete about a basis for $k[\sfX]$ as a $\ph^\ast(k[Y])$-module.  Let $\{x_i \}_1^m$ be coordinate functions for $\sfX$.  By Theorem \ref{thm:geom}, each $x_i$ is integral over $\ph^\ast(k[Y])$.  Let $n_i$ be the degree of an integral dependence for $x_i$ and define $\mathcal{S} = \{ x_1^{\alpha_1} \cdots x_m^{\alpha_m} \, | \, 0\leq \alpha_i < n_i \}$.  Every element of $k[\sfX]$ may be expressed as a $\ph^\ast(k[\sfX])$-linear combination of the elements from the finite set $\mathcal{S}$.  (See Proposition 2.16 of \cite{AM}.) 
\end{enumerate}
\end{remarks}

%%%%%%%%%%%%%%%%%%%%%%%%%%%%%%%%%%%

%%%%%%%%%%%%%%%%%%%%%%%%%%%%%%%%%%%

\section{Character varieties, boundary slopes, and closed essential surfaces} \label{sec:CES}

\begin{dfn}
A {\it knot manifold} is a connected, compact, irreducible, orientable 3-manifold whose boundary is an incompressible torus.
\end{dfn}

Let $N$ be a knot manifold and $\Gamma = \pie(N)$.  Denote the set of $\sl \C$-representations of $\Gamma$ as $\sfR(N)$ and the set of characters of representations in $\sfR(N)$ as $\sfX(N)$.  Let $t \co \sfR(N) \to \sfX(N)$ be the map which takes representations to their characters.  In  \cite{CS1},  it is shown that $\sfR(N)$ and $\sfX(N)$ are affine algebraic sets defined over $\C$ and the map $t$ is regular.  We will refer to $\sfR(N)$ and $\sfX(N)$ as the {\it representation variety} and {\it character variety} for $N$.

Culler and Shalen have revealed deep connections between essential surfaces in $N$ and the character variety $\sfX(N)$.    We outline some of their results in what follows.  For more background, see \cite{CS1}, the survey article \cite{Sh}, or Chapter 1 of \cite{CGLS}.  

\begin{dfns} Suppose that $\sfX$ is a non-empty algebraic subset of $\sfX(N)$.
\begin{enumerate}
\item Given $\gamma \in \Gamma$, the {\it trace function} for $\gamma$ on $\sfX$ is the regular function $I_\gamma \in \C[\sfX]$ defined by $I_\gamma(\chi) = \chi(\gamma)$. 
\item Let $T(\sfX)$ be the subring (with 1) in $\C[\sfX]$ generated by $\{ I_\gamma | \gamma \in \Gamma \}$.  $T(\sfX)$ is called the {\it trace ring} for $\sfX$.  \end{enumerate}
\end{dfns}

The following proposition gives that, as a ring, $T(\sfX)$ is finitely generated.  

\begin{prop}[Proposition 4.4.2 of \cite{Sh}] \label{prop:tracering} Let $\{ \gamma_i \}_1^n$ be a generating set for $\Gamma$.  Then $T(\sfX)$ is generated, as a ring, by the constant function $1$ along with the functions in the set \[ \left\{ I_V  \, \Big| \, V=\gamma_{j_1} \cdots \gamma_{j_k} \text{ where } 1 \leq k \leq n \text{ and } 1 \leq j_1 < \dots < j_k \leq n \right\}.\]
\end{prop}

The smallest $\C$-algebra containing $T(\sfX)$ is the coordinate ring $\C[\sfX]$.  Therefore, a generating set $\{ \gamma_i\}_1^n$ for $\Gamma$ gives an embedding of $\sfX$ into $\C^{2^n-1}$ by taking the functions $I_V$ to be coordinate functions.  It is straightforward to see that under this embedding, $\sfX(N)$ is cut out by polynomials with coefficients in $\Z$.

We have a regular map $\bound \co \sfX(N) \to \sfX(\bound N)$ from $\sfX(N)$ to the character variety for the peripheral subgroup (well-defined up to conjugation) of $\Gamma$ given by restricting characters.  For a non-empty algebraic subset $\sfX$ of $\sfX(N)$, let $\bound \sfX$ denote the Zariski closure of $\bound (\sfX)$.  

\begin{dfns}
\[\] \vspace{-.4in}
\begin{enumerate}
\item The unoriented isotopy class of an essential simple closed curve in $\bound N$ is called a {\it slope}.   We denote the set of slopes on $N$ as $\calS$.
\item A {\it boundary slope} is a slope which is realized as a boundary component of an essential surface in $N$.
\end{enumerate}
\end{dfns}

Each slope corresponds to a pair $\{ \pm a \} \subset H_1(\bound N;\Z)$.  The inverse of the Hurewicz isomorphism is an isomorphism $H_1(\bound N;\Z) \to \pie \bound N$.  This gives a monomorphism $e \co H_1(\bound N;\Z) \to \Gamma$ which is well-defined up to conjugation.  Since traces are invariant under inverses and conjugation, the function $\{ \pm a \} \mapsto I_{e(a)}$ is well-defined.  When $\alpha = \{ \pm a \}$ is a slope we write $I_\alpha = I_{e(a)}$.

Assume now that $\sfX \subseteq \sfX(N)$ is an irreducible affine curve and let \[ \xymatrix{
\wt{\sfX} \ar@{.>}[r]^{\iota \ } &  \sfX^\xi \ar[r]^\xi & \sfX
} \]
be the corresponding maps and varieties as defined in Section \ref{sec:geometry}.   The following theorem is well known and fundamental, see Theorem 2.2.1 and Proposition 2.3.1 of \cite{CS1}.  It is a translation of Theorem \ref{thm:CS} into the language of Sections \ref{sec:geometry} and \ref{sec:CES} of this paper.

\begin{thm}[Culler-Shalen]  \label{thm:cullershalen} For every ideal point $\hat{x}$ of $\sfX$, there is an associated non-empty essential surface $\Sigma$ in $N$.  
\begin{enumerate}
\item $\Sigma$ may be chosen to have empty boundary if and only if  
\[ (\bound \xi \iota)^\ast \left( \C[\bound \sfX] \right) \ \subseteq \ \O_{\wt{\sfX},\hat{x}}.\]
\item Otherwise, there exists a unique slope $\alpha$ such that $I_\alpha \in \O_{\wt{\sfX},\hat{x}}$.  In this case, every component of $\bound \Sigma$  represents the slope $\alpha$.
\end{enumerate}
\end{thm}

Theorem \ref{thm:cullershalen} is broadly applicable.  For instance, if $N$ is the exterior of a non-trivial knot in $S^3$ then Kronheimer and Mrowka show that $\sfX(N)$ contains a curve with infinitely many irreducible characters \cite{KM}.  Also, if the interior of $N$ admits a complete hyperbolic structure with finite volume and $\chi_\rho$ is the character of a discrete faithful representation then there is a unique algebraic component $\sfX_0 \subseteq \sfX(N)$ which contains $\chi_\rho$.  Furthermore, $\sfX_0$ is a curve (see \cite{Th2} and \cite{Sh}).  We refer to such a curve as a {\it hyperbolic curve}.

\begin{dfns}
\[\] \vspace{-.4in}
\begin{enumerate}
\item Suppose that $\sfX \subseteq \sfX(N)$ is an irreducible affine curve and $\hat{x}$ is an ideal point of $\sfX$.  A non-empty essential surface $\Sigma \subset N$ is {\it associated to $\hat{x}$} if it is contained in a surface given by $\hat{x}$ and Theorem \ref{thm:cullershalen}.

\item Suppose that $\sfX$ is an algebraic subset of $\sfX(N)$.  A surface $\Sigma$ is {\it detected by $\sfX$} if there is an ideal point $\hat{x}$ of an irreducible affine curve in $\sfX$ so that $\Sigma$ is associated to $\hat{x}$.

\item If an essential surface $\Sigma \subset N$ is detected by $\sfX$ (associated to $\hat{x}$) and $\bound \Sigma \neq \emptyset$ the boundary slope represented by a component of $\bound \Sigma$ is {\it detected by $\sfX$} (associated to $\hat{x}$).  The slope is {\it strongly} or {\it weakly} detected (associated) depending on whether $\hat{x}$ satisfies case (1) or case (2) of Theorem \ref{thm:cullershalen}, respectively. 
\end{enumerate}
\end{dfns}

\begin{remark}
 It is natural to ask if, whenever $N$ contains a closed essential surface, there is a closed essential surface in $N$ detected by $\sfX(N)$. The author is not certain if the answer is known to be no, however there are compelling reasons to believe that the answer is no, see for example \cite{CT} or Remark 5.1 in \cite{D2}.
\end{remark}

We conclude this section with some applications of the work in Section \ref{sec:geometry}.
\begin{thm} \label{thm:main}
Suppose $N$ is a knot manifold and $\sfX$ is an irreducible algebraic subset of $\sfX(N)$.  The following are equivalent.
\begin{enumerate}
\item $\sfX$ does not detect a closed essential surface.
\item $\dim (\sfX)=\dim (\bound \sfX)=1$ and $\bound \co \sfX \to \bound \sfX$ does not have a hole.
\item $\C[\sfX]$ is integral over $\bound^\ast (\C[\bound \sfX])$.
\item $\C[\sfX]$ is a finitely generated $\bound^\ast(\C[\bound \sfX])$-module.
\end{enumerate}
\end{thm}

\proof As mentioned earlier, conditions (3) and (4) are equivalent by Proposition 5.1 and Corollary 5.2 of \cite{AM}.   

If $\dim(\sfX) > \dim(\bound \sfX)$ then $\C[\sfX]$ is transcendental over $\bound^\ast (\C[\bound \sfX])$ and so condition (3) cannot hold.  The inequality $\dim(\sfX) > \dim (\bound \sfX)$ also implies that $\sfX$ contains a curve $C$ to which Theorem \ref{thm:cullershalen} applies.  This yields  a closed essential surface detected by $\sfX$.  Hence, conditions (2) and (3) both imply that $\dim(\sfX)=\dim(\bound \sfX)$.

We have established that each of the four conditions imply that $\dim(\sfX)=\dim(\bound \sfX)$.  By Proposition 2.4 of \cite{CCGLS}, $\dim(\sfX) \geq 1$.  Also $\dim (\bound \sfX) \leq 1$, otherwise for any slope $\alpha$ on $\bound N$, there is a curve in $\sfX$ to which Theorem \ref{thm:cullershalen} can be applied to give an essential surface in $N$ with boundary slope $\alpha$.  This contradicts  Hatcher's theorem \cite{Hat} that $N$ has only finitely many boundary slopes.  So, if any condition (1) through (4) holds, then $\dim(\sfX)=\dim(\bound \sfX)=1$.

The theorem now follows from Theorems \ref{thm:geom}, \ref{thm:cullershalen}, and Corollary \ref{cor:geom}. \endproof

Given a trace function $I_\gamma \in \C[\sfX]$, let $\C[I_\gamma]$ denote the $\C$-subalgebra of $\C[\sfX]$ generated by $I_\gamma$.  Consider the regular map $I_\gamma \co \sfX \to \C$ and the induced map $I_\gamma^\ast \co \C[x] \to \C[\sfX]$.  If $I_\gamma$ is non-constant on $\sfX$ then $I_\gamma^\ast$ is injective and  $\C[I_\gamma]$ is naturally isomorphic to $\C[x]$.  Moreover, $\Z[x]$ is naturally isomorphic to the $\Z$-submodule $\Z[I_\gamma] \subseteq T(\sfX)$ generated by all finite powers of $I_\gamma$.  For instance, if $\alpha$ is a slope which is not detected by $\sfX$, since $\dim(\sfX) \geq 1$, Theorem \ref{thm:cullershalen} implies that $I_\alpha$ is non-constant on $\sfX$.

\begin{thm} \label{thm:I_mu}
Suppose $N$ is a knot manifold, $\sfX$ is an irreducible algebraic subset of $\sfX(N)$, and $\alpha$ is a slope.  The following are equivalent.
\begin{enumerate}
\item $\sfX$ does not detect a closed essential surface and $\alpha$ is not detected by $\sfX$.
\item $\C[\sfX]$ is integral over $\C[I_\alpha]$.
\item As a $\C[I_\alpha]$-module, $\C[\sfX]$ is a finitely generated and free.
\end{enumerate}
\end{thm}

\proof As with conditions (3) and (4) of the previous theorem, we know that (3) implies (2) and that (2) implies that $\C[\sfX]$ is a finitely generated $\C[I_\alpha]$-module.  $\C[\sfX]$ is an integral domain so, as a $\C[I_\alpha]$-module, $\C[\sfX]$ is torsion free.  $\C[I_\alpha]$ is a PID so $\C[X]$ is a free $\C[I_\alpha]$-module.

It remains to show that (1) and (2) are equivalent.  

Assume first that (2) holds.  Theorem \ref{thm:main} shows that $\sfX$ does not detect a closed essential surface and $\dim (\sfX)=1$.  Hence $I_\alpha$ is not constant.  We have the regular map $I_\alpha \co \sfX \to \C$ and $\overline{I_\alpha(\sfX)}=\C$.  Since  $\C[\sfX]$ is integral over $\C[I_\alpha]$, Theorem \ref{thm:geom} shows that $I_\alpha \co \sfX \to \C$ does not have a hole.  By Theorem \ref{thm:cullershalen}, $\alpha$ is not strongly detected by $\sfX$.  Since $\sfX$ does not detect a closed essential surface, $\alpha$ cannot be weakly detected by $\sfX$.

Now suppose that (1) holds.  As in the proof of Theorem \ref{thm:main}, we have that $\dim(\sfX)=1$.  Consider the regular map $I_\alpha \co \sfX \to \C$.  Either $\overline{I_\alpha(\sfX)}=\C$ or $I_\alpha$ is constant on $\sfX$.   But since $\dim(\sfX)=1$, Theorem \ref{thm:cullershalen} implies that if $I_\alpha$ is constant on $\sfX$ then $\alpha$ is strongly detected by $\sfX$ or $\sfX$ detects a closed essential surface.  So we must have $\overline{I_\alpha(\sfX)}=\C$.

By Theorem \ref{thm:geom}, it suffices to show that the regular map $I_\alpha$ does not have a hole.  Again, we appeal to Theorem \ref{thm:cullershalen} and notice that if $I_\alpha$ has a hole then either $\alpha$ is strongly detected or $\sfX$ detects a closed essential surface. 
\endproof

\begin{dfn} \label{dfn:Crank} Let $\sfX \subseteq \sfX(N)$ be an irreducible algebraic subset.    The $\C$-{\it rank function} is the map \[\rk_\sfX^\C \co \calS \to \Z^+ \cup \{ \infty\}\] where $\rk_\sfX^\C(\alpha)$ is the rank of $\C[\sfX]$ as a $\C[I_\alpha]$-module. \end{dfn}

\begin{remarks} Using Theorem \ref{thm:I_mu}, we make the following observations.
%\[\] \vspace{-.4in}
\begin{enumerate}
\item If $\sfX$ does not detect a closed essential surface, $\rk_\sfX^\C(\alpha)=\infty$ if and only if $\alpha$ is strongly detected by $\sfX$.
\item $\sfX$ detects a closed essential surface if and only if $\tr_\sfX^\C(\calS)=\{\infty\}$.  
\item By Lemma 1.4.4 of \cite{CGLS}, there are only finitely many boundary slopes strongly detected by $\sfX$.  Avoiding this finite set, choose $\alpha \in \calS$.  Then $\sfX$ detects a closed essential surface if and only if $\rk_\sfX^\C(\alpha) = \infty$.
\end{enumerate}
\end{remarks}

For $\sfX$ an irreducible algebraic component of $\sfX(N)$ we define $\sfX^\circ$ to be the union of the irreducible components $\sfX'$ of $\sfX(N)$ such that
\begin{enumerate}
\item $\sfX'$ detects a closed essential surface if and only if $\sfX$ detects a closed essential surface and 
\item for every slope $\alpha$, $\sfX'$ detects $\alpha$ if and only $\sfX$ detects $\alpha$.
\end{enumerate}

\begin{prop} \label{prop:X^circ}
If $\sfX$ is an irreducible algebraic component of $\sfX(N)$ 
then $\sfX^\circ$ is defined over $\Q$.  
\end{prop}
\proof
By a theorem of Weil \cite[Ch. III, Thm. 7]{La}, it suffices to prove that $\sfX^\circ$ is invariant under the action of $\text{Aut}(\C)$ on $\sfX(N)$ (see \cite[Prop. 2.3]{CCGLS} and \cite[Section 5]{BZ2}).  

First, we assume that $\sfX$ does not detect a closed essential surface.   Theorem \ref{thm:I_mu} implies that $\sfX^\circ$ is the union of the irreducible algebraic components $\sfX'$ of $\sfX(N)$ such that, for every slope $\alpha$, $\C[\sfX']$ is integral over $\C[I_\alpha]$ if and only if $\C[\sfX]$ is integral over $\C[I_\alpha]$. 

Let $\phi \in \text{Aut}(\C)$, $\sfX'$ an irreducible component of $\sfX^\circ$, and $\alpha \in \calS$.  Define $\sfX''=\phi(\sfX')$.  The set $\phi(\sfX')$ is an irreducible algebraic component of $\sfX(N)$. To prove the proposition we need only show that $\C[\sfX'']$ is integral over $\C[I_\alpha]$ if and only if $\C[\sfX']$ is integral over $\C[I_\alpha]$.

The automorphism $\phi$ determines an automorphism on any complex polynomial ring by acting on the coefficients.  This automorphism descends to an isomorphism $\phi \co \C[\sfX'] \to \C[\sfX'']$. Since $\phi \co \C \to \C$ restricts to the identity on $\Q$ and $I_\alpha$ is represented by a polynomial with integer coefficients we have $\phi(I_\alpha)=I_\alpha$.  Therefore if $p \in \C[I_\alpha]$ then $\phi(p) \in \C[I_\alpha]$.  So if $f \in \C[\sfX']$ is integral over $\C[I_\alpha]$, we can apply $\phi$ to an integral dependence relation for $f$ over $\C[I_\alpha]$ to obtain one for $\phi(f)$  over $\C[I_\alpha]$.  To prove the converse, we apply $\phi^{-1}$ to any integral dependence relation for $\C[\sfX'']$ over $\C[I_\alpha]$.  

The case when $\sfX$ does detect a closed essential surface follows from the above argument along with Lemma 5.3 of \cite{BZ2}.
\endproof

\begin{dfn} Suppose $\sfX \subseteq \sfX(N)$ is an algebraic set.  The {\it rational trace ring} $T_\Q(\sfX)$ is the smallest $\Q$-algebra that contains $T(\sfX)$.  
\end{dfn}

\begin{thm} \label{thm:Q}
Suppose $N$ is a knot manifold, $\sfX$ is an irreducible algebraic subset of $\sfX(N)$, and $\alpha \in \calS$.  The following are equivalent.
\begin{enumerate}
\item $\sfX$ does not detect a closed essential surface and $\alpha$ is not detected by $\sfX$.
\item $T_\Q(\sfX)$ is integral over $\Q[I_\alpha]$.
\item As a $\Q[I_\alpha]$-module, $T_\Q(\sfX)$ is a finitely generated and free.
\end{enumerate}
\end{thm}

\proof
The first paragraph of the proof of Theorem \ref{thm:I_mu} shows that (2) and (3) are equivalent.  Theorem \ref{thm:I_mu} shows that (2) implies (1).  

It remains to show that (1) implies (2).  Suppose then that $\sfX$ does not detect a closed essential surface and $\alpha$ is not detected by $\sfX$.  $T_\Q(\sfX)$ is generated by $T(\sfX)$ as a $\Q$-algebra, so it suffices to show that elements of $T(\sfX)$ are integral over $\Q[I_\alpha]$. 

Let $\{ \sfX_j \}_1^n$ be the set of irreducible algebraic components of $\sfX^\circ$ and take $f \in T(\sfX)$.  The map $\C[\sfX^\circ] \to \C[\sfX]$ induced by inclusion is given by restriction and is surjective.  Take $F \in \C[\sfX^\circ]$ such that $F|_\sfX=f$.  Let $\pi \co \sfX^\circ \to \C^2$ be the map $(F,I_\alpha)$.  Let $(x,y)$ be the coordinates on $\text{Im}(\pi)$ determined by $x\pi=F$ and $y\pi=I_\alpha$.

Fix $j \in \{ 1, \ldots , n\}$.  By Theorem \ref{thm:main},  $I_\alpha \co \sfX_j \to \C$ is non-constant and $\overline{\pi(\sfX_j)}$ is an irreducible plane curve.  Let $P_j \in \C[x,y]$, be an irreducible polynomial which defines $\overline{\pi(\sfX_j)}$.   $\sfX^\circ$ is defined over $\Q$ (Proposition \ref{prop:X^circ}) and $\pi$ is given by polynomials with $\Z$-coefficients, so we may assume that \[ P \ = \ \prod_{j=1}^n P_j \ \in \Z[x,y].\]  

For each $j$, Theorem \ref{thm:I_mu} gives that $F|_{\sfX_j}$ is integral over $\C[I_\alpha]$.  This means that there is a polynomial in $\C[x,y]$ which is zero on $\overline{\pi(\sfX_j)}$ and monic as a polynomial in $(\C[x])[y]$.   Select one such polynomial $p_j$ with minimal degree. Then $p_j$ is irreducible and must differ from $P_j$ only by multiplication by some $\alpha_j \in \C$.  We have \[ P \ = \ \prod P_j \ = \ \left(\prod \alpha_j \right) \cdot \left( \prod p_j\right).\] Hence $\prod \alpha_j \in \Z$ and $(\prod \alpha_j^{-1}) \cdot P$ gives an integral dependence relation for $f$ over $\Q[I_\alpha]$.
\endproof

\begin{dfn} \label{dfn:Qrank} Suppose $\sfX \subseteq \sfX(N)$ is an irreducible algebraic subset.    The $\Q$-{\it rank function} is the map \[\rk_\sfX^\Q \co \calS \to \Z^+ \cup \{ \infty\}\] where $\rk_\sfX^\Q(\alpha)$ is the rank of $T_\Q(\sfX)$ as a $\Q[I_\alpha]$-module. \end{dfn}

It will require more work to show that there are related theorems which work over $\Z$.  We begin by establishing some lemmas concerning valuations.

%%%%%%%%%%%%%%%%%%%%%%%%%%%%%%%%%%%

%%%%%%%%%%%%%%%%%%%%%%%%%%%%%%%%%%%

\section{Valuations}

\begin{lemma} \label{lem:val2}
Suppose $\nu \co F^\times \to \Lambda$ is a valuation on the field $F$.  Assume that \[ cx^k \ = \ \sum_{i=0}^{k-1} a_i x^i\]
for some $c,x,a_i \in F$ and $k\in \Z^+$ such that $c \neq 0$, $\nu(c)=0$, and $\nu(a_i)\geq0$ for every $i$.  Then $\nu(x) \geq 0$.
\end{lemma}
\proof
If every $a_i$ is zero then $x=0$ and so $\nu(x)=\infty$, so we may assume that at least one of the $a_i$'s is non-zero.  Let $j$ be such that $\nu(a_j)+j\nu(x)$ is minimal in the set $\{ \nu(a_i)+i\nu(x) \, | \, 0\leq i <k \}$.  We have
\[ k \nu(x) \ = \ \nu(c)+k\nu(x) \ \geq \ \nu(a_j)+j\nu(x) \]
so
\[ (k-j)\nu(x)\  \geq \ \nu(a_j)\ \geq \ 0 .\]
Hence $\nu(x) \geq 0$.
\endproof

\begin{lemma} \label{lem:padic}
If $\nu$ is a non-trivial valuation on $\Q$ then $\nu$ is a $p$-adic valuation.
\end{lemma}
\proof
Let $R$ be the valuation ring corresponding to $\nu$ and let $\mathfrak{m}$ be the maximal ideal in $R$.  We know that $\Z \subset R$.  Note that $\mathfrak{m} \cap \Z$ is prime (maximal) in $\Z$ so $\nu$ is $p$-adic. 
\endproof

The following lemma is essentially Lemma 2.2 of \cite{CL2}.
\begin{lemma} \label{lem:CL}
Suppose $p \in \Z$ is a prime and $h\in \Z[t]$ is irreducible over $\Z$ and of the form
 \[ h(t) \ = \ c_0 + c_1 t+ \dots +c_{k-1}t^{k-1}+c_kp^rt^k \] where $(c_k,p)=1$ and $c_k \neq \pm 1$.  Let $b\in \C$ be any root of $h$.  Then there is a valuation $\nu$ on $\Q(b)$ with $\nu(p)=0$ and $\nu(b)<0$.
\end{lemma}

\proof
Since $(c_k,p)=1$ and $h\in \Z[t]$ we know that $r$ is an integer and at least zero.  Now, for $i \in \{0, \ldots, k-1\}$, define $a_i=c_ip^{k-1-i}$.  We have $p^{(k-1)r}c_i=a_i(p^r)^i$, so
\[p^{(k-1)r} h(t) \ =\ a_0+a_1(p^rt)+ \dots +a_{k-1}(p^rt)^{k-1}+c_k(p^rt)^k \]
Set $d=p^rb$ and \[ h^\ast(t) \ =\  a_0  + a_1 t+ \dots +a_{k-1}t^{k-1}+c_kt^k.\]  Then $h^\ast \in \Z[t]$ and $h^\ast(d)=0$.   We have
\[ c_kd^k \ = \ \sum_{i=0}^{k-1} a_i d^i.\]
Let $\nu_p$ be an extension of the $p$-adic valuation to $\Q(b)$.  Since $(c_k,p)=1$, $\nu_p(c_k)=0$.   Also, for every $i$, $a_i \in \Z$ so $\nu(a_i)\geq 0$.  Hence, Lemma \ref{lem:val2} implies that $\nu_p(d)\geq0$.

Since $h$ is irreducible over $\Z$ and $\Q(b)=\Q(d)$, it follows that $h^\ast$ is irreducible over $\Z$.  Hence, neither $b$ nor $d$ is integral over $\Z$.  By Corollary 5.2 of \cite{AM}, we have a valuation $\nu$ on $\Q(b)$ with $\nu(d)<0$.  The restriction of $\nu$ to $\Q$ cannot be the $p$-adic valuation and (again using Lemma \ref{lem:val2}) cannot be trivial.   By Lemma \ref{lem:padic}, $\nu|_{\Q}=\nu_q$ for some prime $q \in \Z$ which is different from $p$.  Therefore, $\nu(p)=0$.   Also
\[ \nu(b) \ = \ \nu(dp^{-r}) \ = \ \nu(d)-r\nu(p)\ =\ \nu(d) <0.\]
\endproof

%%%%%%%%%%%%%%%%%%%%%%%%%%%%%%%%%%%

%%%%%%%%%%%%%%%%%%%%%%%%%%%%%%%%%%%

\section{$T(\sfX)$ as a $\Z[I_\alpha]$-module} \label{sec:Z}

Throughout this section, $N$ is a knot manifold and $\Gamma$ is its fundamental group.  Next, we introduce some terminology from \cite{SZ}.  

\begin{dfn}
A representation $\rho \co \Gamma \to \sl \C$ is {\it algebraic non-integral}  if the image of $\rho$ is in $\sl F$, where $F$ is a number field and $\chi_\rho(\gamma)$ is not an algebraic integer for some $\gamma \in \Gamma$.  If $\rho$ is an algebraic non-integral representation we abbreviate this by saying that $\rho$ is an {\it ANI-representation}.

\end{dfn}

Similar to Theorem \ref{thm:cullershalen}, the following result again follows from the work of Culler and Shalen.  See Lemma 2 of \cite{SZ}.

\begin{thm}  \label{thm:cullershalen_ANI} For every ANI-representation $\rho$ of $\Gamma$ there is an associated non-empty essential surface $\Sigma \subset N$.  Suppose that $\alpha \in \calS$ is not a boundary slope.  The surface $\Sigma$ may chosen to have empty boundary if and only if  $I_\alpha(\chi_\rho)$ is an algebraic integer.
\end{thm}

\begin{dfns}
\[\] \vspace{-.4in}
\begin{enumerate}
\item Suppose that $\rho$ is an ANI-representation of $\Gamma$ and $\Sigma \subset N$ is a non-empty essential surface.  $\Sigma$ is {\it associated to $\rho$} if it is contained in a surface given by $\rho$ and Theorem \ref{thm:cullershalen_ANI}.   

\item Suppose that $\sfX$ is an algebraic subset of $\sfX(N)$.  A surface $\Sigma$ is {\it ANI-detected by $\mathbf{\sfX}$} if there is an ANI-representation $\rho$ of $\Gamma$ whose character lies on $\sfX$ and  $\Sigma$ is associated to $\rho$.
\end{enumerate}
\end{dfns}
 
\begin{prop} \label{prop:ANI}
If $\sfX$ is an irreducible algebraic subset of $\sfX(N)$ and $N$ contains a closed essential surface which is ANI-detected by $\sfX$, then $N$ contains a closed essential surface detected by $\sfX$.
\end{prop}

\proof
Assume that $\sfX$ does not detect a closed essential surface and suppose for a contradiction that there is a closed essential surface in $N$ associated to an ANI-representation $\rho \co \Gamma \to \sl F$ with $\chi_\rho \in \sfX$.  Take $\gamma \in \Gamma$ with $\chi_\rho(\gamma)$ non-integral and choose $\alpha \in \calS$ not a boundary slope.  

By Theorems \ref{thm:intersect} and \ref{thm:cullershalen_ANI}, there is a valuation $\nu$ on $F$ with $\nu(I_\gamma(\chi_\rho))<0$ and $\nu(I_\alpha(\chi_\rho))\geq 0$.  By Theorem \ref{thm:Q}, there is a polynomial \[ f(x,y) \ = \ c x^r + \sum_{j=0}^{r-1} q_j(y)x^j \] with $c \in \Z-\{0\}$, $q_j \in \Z[y]$, and $f(I_\gamma, I_\alpha)=0$ on $\sfX$.  Since $\nu(I_\alpha(\chi_\rho))\geq 0$, we know that $\nu(q_j(I_\alpha(\chi_\rho)))\geq 0$.  But then, Lemma \ref{lem:val2} implies that $\nu(I_\gamma(\chi_\rho)) \geq 0$, a contradiction.
\endproof

\begin{example}  In \cite{CL3} Cooper and Long ask whether there is a one-cusped manifold whose character variety detects a closed essential surface.  This was answered by Tillmann in \cite{Ti} where he presented a calculation which shows that the character variety of the Kinoshita-Terasaka knot (11n42) detects a closed essential surface.  Together with data collected by Goodman, Heard, and Hodgson  \cite{GHH}, \url{ http://www.ms.unimelb.edu.au/~snap/  }, Proposition \ref{prop:ANI} gives many more examples. 

Suppose that $N$ is a one-cusped hyperbolic 3-manifold and let $\sfX_0$ be an algebraic component of $\sfX(N)$ which contains the character of a discrete faithful representation $\rho_0$.  By Theorem \ref{thm:cullershalen_ANI} and Proposition \ref{prop:ANI}, if $\rho_0$ is an ANI-representation then $N$ contains a closed essential surface detected by $\sfX(N)$.  

For all manifolds in the Callahan-Hildebrand-Weeks census of cusped hyperbolic manifolds with up to $7$ tetrahedra \cite{CHW}, and for all complements of hyperbolic knots and links up to 12 crossings, Goodman, Heard, and Hodgson use the computer program {\tt Snap} to attempt to (among other things) decide whether or not the holonomy representation is an ANI-representation.  They find 252 such manifolds, 21 of which are complements of knots in $S^3$ (this is almost certainly not a complete list).   The knots are 9a30,  9a31,  10a89,   10a96,  10a103, 11n97, 12n156, 12n245, 12n246, 12n260, 12n494, 12n508, 12n518, 12n600, 12n602, 12n604, 12n605, 12n694, 12n888, 12a1205, and 12a1288.
\end{example}

Roughly, the following theorem is proven by applying Theorem \ref{thm:Q} and arguing as in \cite{CL2}.

\begin{thm} \label{thm:Z}
Suppose $N$ is a knot manifold, $\sfX$ is an irreducible algebraic subset of $\sfX(N)$, and $\alpha \in \calS$.  The following are equivalent.
\begin{enumerate}
\item $\sfX$ does not detect a closed essential surface and $\alpha$ is not detected by $\sfX$.
\item $T(\sfX)$ is integral over $\Z[I_\alpha]$.
\item $T(\sfX)$ is a finitely generated $\Z[I_\alpha]$-module.
\end{enumerate}
\end{thm}

\proof
As with Theorem \ref{thm:Q}, we need only argue that (1) implies (2).  Suppose that (1) is true and take $\gamma \in \Gamma$.  Our goal is to show that $I_\gamma$ is integral over $\Z[I_\alpha]$, so we may assume that $I_\gamma$ is not everywhere zero on $\sfX$.  By Theorem \ref{thm:main}, $\dim(\sfX)=1$ and by Theorem \ref{thm:I_mu}, $I_\alpha$ must be non-constant on $\sfX$.

By Theorem \ref{thm:Q}, there is a polynomial $f(y,z) \in \Z[y,z]$ of the form \[ f(y,z) \ =\ c y^r + \sum_{j=0}^{r-1} q_j(z)y^j \]
with $q_j \in \Z[z]$ such that $c \neq 0$, $f(I_\gamma, I_\alpha)$ is the zero function on $\sfX$, and the greatest common divisor of the coefficients of $f$ is one.  Choose $f$ among all such polynomials to have smallest possible total degree.  Then, since $I_\alpha$ and $I_\gamma$ are not everywhere zero on $\sfX$, neither $y$ nor $z$ divides $f$.

Take $k \in \mathbb{N}$ large enough such that $f(y,xy^{-k})$ has no two terms with the same degree in $y$.  Take $l \in \Z$ as small as possible such that \[ g(x,y) \ = \ y^l f(y,xy^{-k}) \] is a polynomial in $\Z[x,y]$.  Then $g(x,y)$ is of the form \[ g(x,y) \ =\ c y^s + \sum_{j=0}^{s-1} a_j x^{t_j} y^j \] where the $t_j$'s are distinct integers, the greatest common divisor of the $a_j$'s is one, $g(I_\alpha I_\gamma^k, I_\gamma)=0$ on $\sfX$, and neither $x$ nor $y$ divides $g$.

We know that $I_\alpha$ is non-constant on $\sfX$ so there must be a point $\hat{x} \in \wt{\sfX}$ with $I_\alpha(\hat{x})=0$.  If $I_\alpha I_\gamma^k$ is constant on $\sfX$ then $I_\alpha= K I_\gamma^{-k}$ for some $K \in \C$.   Hence, $\hat{x}$ is an ideal point of $\sfX$.
 By Theorem 2.2.1 and Proposition 2.3.1 of \cite{CS1}, either $\sfX$ detects a closed essential surface or an essential surface with boundary slope $\alpha$.  Both are contradictions so $I_\alpha I_\gamma^k$ must be non-constant on $\sfX$.  Therefore, we may choose a prime $p$ such that $\gcd(c , p)=1$ and for every root $b$ of $g(p,y)$ there is a character $\chi_\rho \in \sfX$ with $\chi_\rho(\alpha)\cdot \chi_\rho(\gamma)^k = p$ and $\chi_\rho(\gamma)=b$.
 
 Take \[ g(p,y) \ = \ \eta \cdot \prod_{j=0}^u h_j(y) \] a complete factorization of $g(p,y)$ over $\Z$.  So $\eta \in \Z$ and each $h_j$ is irreducible over $\Z$.  We claim that $c$ divides $\eta$.  The leading coefficient of $g(p,y)$ is $c$.  If we denote the the leading coefficient of $h_j$ as $c_j$, we have $c = \eta c_0 \cdots c_u$.  If $c$ doesn't divide $\eta$ then there is some $c_i$ with $\gcd(\alpha, c_i)>1$.  Write $c_i = p^t d$ where $\gcd(p,d)=1$.  Since $\gcd(c , p)=1$ we must have $\gcd(c, d)>1$, in particular $d \neq \pm 1$.   Let $b$ be a root of $h_i$.  By our choice of $p$, there is $\chi_\rho \in \sfX$ which corresponds to the solution $g(p,b)$.  Using Hilbert's Nullstellensatz, we may assume that the entries of the matrices in $\rho(\pie N)$  lie in a number field $F$.  By Lemma \ref{lem:CL}, we have a valuation $\nu$ on $\Q(b)$ with $\nu(p)=0$ and $\nu(b)<0$.  Take an extension of $\nu$ to $F$.  We have $\nu(I_\gamma(\chi_\rho))<0$ and \[ 0 \ =\ \nu\left(I_\alpha(\chi_\rho) I_\gamma(\chi_\rho)^k\right) \ =\ \nu\left(I_\alpha(\chi_\rho) \right) + k \cdot \nu\left( I_\gamma(\chi_\rho) \right).\]   Since $k \geq 0$ we must have $\nu(I_\alpha(\chi_\rho))\geq 0$.  Hence, $N$ contains a closed essential surface which is ANI-detected by $\sfX$.  Proposition \ref{prop:ANI} shows that this is a contradiction, so we must have that $c$ divides $\eta$.
 
 The constant $\eta$ divides every coefficient of $g(p,y)$ and $c$ divides $\eta$, so $c$ divides every coefficient of $g(p,y)$.  But the coefficient for $z^j$ is $a_j p^{t_j}$ and $\gcd(c, p)=1$.  Hence $c$ divides each $a_j$.  The only integers that divide every $a_j$ are $\pm 1$.  Therefore $f(I_\gamma, I_\alpha)$ is an integral dependence relation for $I_\gamma$ over $\Z[I_\alpha]$.
\endproof

\begin{dfn} \label{dfn:Zrank} Let $\sfX \subseteq \sfX(N)$ be an irreducible algebraic subset.    The $\Z$-{\it rank function} is the map \[\rk_\sfX^\Z \co \calS \to \Z^+ \cup \{ \infty\}\] where $\rk_\sfX^\Z(\alpha)$ is the rank of $T(\sfX)$ as a $\Z[I_\alpha]$-module. \end{dfn}

\begin{cor}\label{cor:irred}
Suppose $N$ is a knot manifold, $\sfX$ is an irreducible algebraic subset of $\sfX(N)$ which is defined over $\Q$.  The following are equivalent.
\begin{enumerate}
\item $\sfX$ does not detect a closed essential surface and $\alpha \in \calS$ is not detected by $\sfX$.
\item Every $f \in T(\sfX)$ satisfies a monic \underline{irreducible} polynomial with coefficients in $\Z[I_\alpha]$. 
\end{enumerate}
\end{cor}

\proof 
By Theorem \ref{thm:I_mu}, need only argue that (1) implies (2).  Assume that $\sfX$ does not detect a closed essential surface and let $f$ be a non-zero element of  $T(\sfX)$.  By Theorem \ref{thm:main}, $\sfX$ is a curve. Since $\sfX$ does not detect a closed essential surface, $I_\alpha$ is non-constant on $\sfX$.  Let $V$ be the Zariski closure of the image of $\sfX$ under the map $(f,I_\alpha)$.  $V$ must be an irreducible affine curve defined over $\Q$.  Let $p \in \Q[x,y]$ be a defining equation for $V$.  By multiplying by an integer we may assume $p \in \Z[x,y]$.  The proof of Theorem \ref{thm:Z}, shows that $p$ is an integral dependence relation of $f$ over $\Z[I_\alpha]$. Since $V$ is irreducible, $p$ is irreducible.
\endproof

\begin{dfn} A knot manifold is called {\it small} if it does not contain a closed essential surface.  \end{dfn}

\begin{cor} \label{cor:small}
Let $\sfX$ be an irreducible algebraic subset of $\sfX(N)$, where $N$ is a small knot manifold.  A slope $\alpha$ is detected by $\sfX$ if and only if $\rk_\sfX^\Z(\alpha) \neq \infty$.  
\end{cor}

%%%%%%%%%%%%%%%%%%%%%%%%%%%%%%%%%%%

%%%%%%%%%%%%%%%%%%%%%%%%%%%%%%%%%%%

\section{$\psl \C$-character varieties}

The theorems from Sections \ref{sec:CES} and \ref{sec:Z}, also apply in the $\psl \C$ setting after making appropriate definitions and using virtually identical arguments.  See \cite{BZ1} for details in what follows.

As usual, let $N$ be a compact irreducible orientable 3-manifold with torus boundary and let $\Gamma=\pie N$.  Denote the set of $\psl \C$-representations of $\Gamma$ as $\overline{\sfR}(N)$.  $\psl \C$ acts on $\overline{\sfR}(N)$ by conjugation.  Let $\bar{t} \co \overline{\sfR}(N) \to \sfY(N)$ denote the corresponding algebro-geometric quotient.  $\overline{\sfR}(N)$ and $\sfY(N)$ are affine algebraic sets and the map $\bar{t}$ is a surjective regular map.  $\sfY(N)$ is called the $\psl \C$-character variety for $N$ and a point in $\sfY(N)$ is called a $\psl \C$-character.  As with the $\sl \C$-character variety, every irreducible algebraic component of $\sfY(N)$ has dimension at least one \cite{CCGLS}.

%The natural quotient $\sl \C \to \psl \C$ induces a regular map $\Phi \co \sfX(N) \to \sfY(N)$.  For $\bar{\chi} \in \sfY(N)$, $\Phi^{-1}(\bar{\chi})$ is either empty or is the orbit of a point in $\sfX(N)$ under the action of $H^1(N, \Z_2)$ on $\sfY(N)$.  Here we view $\Z_2$ as the multiplicative group $\{ \pm 1\}$ and for $\sigma \in H^1(N;\Z_2)$ and $\gamma \in \Gamma$, $\sigma$ acts by  $(\sigma \cdot \chi)(\gamma) = \sigma(\gamma) \chi(\gamma)$.

We would like to relate $\sfY(N)$ to $\sl \C$-character varieties and define $T(\sfY)$ where $\sfY$ is an irreducible algebraic subset of $\sfY(N)$.  First we claim that there is a well-defined equivalence class of $\Z_2$-central extensions \begin{equation} \label{eq:extension} 1 \to \Z_2 \to \hat{\Gamma} \to \Gamma \to 1 \end{equation}
and a finite regular map from the $\sl \C$-character variety $\sfX(\hat{\Gamma})$ of $\hat{\Gamma}$ to $\sfY(N)$ which contains $\sfY$ in its image.

First note that if $\bar{\rho}_1, \bar{\rho_2} \in \overline{\sfR}(N)$ have the same image under $\bar{t}$ then the representations are either both reducible or are both irreducible.   Hence we may refer to the points of $\sfY$ as being either reducible or irreducible.   

Suppose that every element of $\sfY$ is reducible and choose $\bar{\chi} \in \sfY$ smooth in $\sfY(N)$.  There is a diagonal representation $\bar{\rho} \in \bar{t}^{-1}(\bar{\chi})$.  This representation determines an extension (\ref{eq:extension}) and a lift $\hat{\rho} \co \hat{\Gamma} \to \sl \C$ of $\bar{\rho}$.  Moreover, the isomorphism class of $\hat{\Gamma}$ is independent of our choices of $\bar{\chi}$ and $\bar{\rho}$.  The natural epimorphisms $\sl \C \to \psl \C$ and $\hat{\Gamma} \to \Gamma$ together induce a regular map $\phi \co \sfX(\hat{\Gamma}) \to \sfY(N)$.  In fact, the image of $\phi$ is the quotient of $\sfX(\hat{\Gamma})$ under the natural action of $H^1(\hat{\Gamma};\Z_2)$.  Moreover, $\sfY \subseteq \text{Im}(\phi)$.

Otherwise $\sfY$ contains an irreducible character which is a smooth point of $\sfY(N)$.  Let $\bar{\chi}$ be such a point and $\bar{\rho} \in \bar{t}^{-1}$.  Exactly as before, we obtain an extension (\ref{eq:extension}) and a lift $\hat{\rho}$ of $\bar{\rho}$ to $\hat{\Gamma}$.  Again, the isomorphism class of $\hat{\Gamma}$ is independent of our choices, we obtain a finite regular map $\phi \co \sfX(\hat{\Gamma}) \to \sfY(N)$, the image of $\phi$ is the quotient of $\sfX(\hat{\Gamma})$ by $H^1(\hat{\Gamma};\Z_2)$, and  $\sfY\subseteq \text{Im}(\phi)$.

In either case, the induced map $\phi^\ast \co \C[\text{Im}(\phi)] \to \C[\sfX(\hat{\Gamma})]$ is injective.  We define the trace ring $T(\text{Im}(\phi))$ to be \[ T( \text{Im}(\phi) ) \ = \  (\phi^\ast)^{-1} \left( \text{Im}(\phi^\ast) \cap T\big(\sfX(\hat{\Gamma})\big) \right).\]  The inclusion $\sfY \to \text{Im}(\phi)$ induces an epimorphism $ \C[\text{Im}(\phi)] \to  \C[\sfY]$.  Define $T(\sfY)$ to be the image of $T( \text{Im}(\phi) )$ under this epimorphism.  

Note that if $\gamma \in \Gamma$ and $\hat{\gamma}_1, \hat{\gamma_2} \in \hat{\Gamma}$ are the preimages of $\gamma$ then $I_{\hat{\gamma}_1}^2=I_{\hat{\gamma}_2}^2$.  It follows that there is a well-defined squared trace function $I_\gamma^2 \in T(\sfY)$.

To obtain results in the $\psl \C$-setting corresponding to those in  Sections \ref{sec:CES} and \ref{sec:Z}, we use the $\psl \C$-definitions rather than $\sl \C$-definitions and replace $I_\alpha$ with the squared trace function $I_\alpha^2$ whenever $I_\alpha$ is not in the $\psl \C$-trace ring.  

\begin{remark} The equivalence classes of extensions (\ref{eq:extension}) are in bijective correspondence with the elements of $H^2(\Gamma; \Z_2)$.
\end{remark}

%The equivalence classes of $\Z_2$-central extensions of $\Gamma$ correspond to the elements of $H^2(\Gamma; \Z_2)$.  Let $\sigma$ be a 2-cocycle which represents the class of (\ref{eq:extension}).  Then $\hat{\Gamma}$ is isomorphic to the group $E_\sigma = \Z_2 \times \Gamma$, where \[ (a,g)(b,h)= \left(a+b+\sigma(g,h), gh\right) \AND (a,g)^{-1}=\left(a-\sigma(g^{-1},g),g^{-1}\right).\]  It is straightforward to check that if $\sigma'$ is another 2-cocycle defining an extension $E_{\sigma'}$ with a lift $\rho' \co E_{\sigma'} \to \sl \C$ of $\bar{\rho}$ then $\sigma-\sigma'$ is zero in $H^2(\Gamma;\Z_2)$ and that $E_\sigma \cong E_{\sigma'}$.  Hence if $\rho'$ is any lift of $\bar{\rho}$ to a $\Z_2$-central extension then $\chi_{\rho'} \in \sfX(\hat{\Gamma})$.  
%
%The quotient $\sl \C \to \psl \C$ induces a surjective regular map $\sfX(\hat{\Gamma}) \to \sfY(\hat{\Gamma})$.   The epimorphism $\hat{\Gamma} \to \Gamma$ induces a regular map $\sfY(\hat{\Gamma}) \to  \sfY(N)$.  Composition gives a regular map $\phi \co \sfX(\hat{\Gamma}) \to \sfY(N)$.  This map is the quotient of $\sfX(\hat{\Gamma})$ by $H^1(\hat{\Gamma};\Z_2)$ and $\sfY \subseteq \text{Im}(\phi)$.  The induced map $\phi^\ast \co \C[\text{Im}(\phi)] \to \C[\sfX(\hat{\Gamma})]$ is injective.  

%%%%%%%%%%%%%%%%%%%%%%%%%%%%%%%%%%%

%%%%%%%%%%%%%%%%%%%%%%%%%%%%%%%%%%%

\section{Ranks and examples}

Throughout this section we assume that $N$ is a knot manifold and $\Gamma=\pie N$.  For an irreducible algebraic set $\sfX \subseteq \sfX(N)$, observe that \begin{equation} \label{eq:inequality} \rk_\sfX^\C \ \leq \ \rk_\sfX^\Q \ \leq \ \rk_\sfX^\Z.\end{equation} Among other things, this section gives examples which show that these inequalities need not be equalities.

If $R$ is a ring, $M$ is an $R$-module, and $S \subset M$ then we write $\langle S \rangle_R$ to denote the $R$-module generated by $S$.

%%%%%%%%%%%%%%%%%%%%%%%%%%%%%%%%%%%

\subsection{Abelian curves}  Suppose that $H_1(N;\Z) \cong \Z$.  There is a curve $\sfXA \subseteq \sfX(N)$ consisting of all characters of representations which factor through $H_1(N;\Z)$.  Let $\alpha \in \Gamma$ be a peripheral element whose image in $H_1(N;\Z)$ generates $H_1(N;\Z)$.  

Every character in $\sfXA$ is the character of a diagonal representation and the image subgroup in $\sl \C$ is generated by the image of $\alpha$.  It follows that every trace function $I_\alpha \in T(\sfXA)$ is represented by a polynomial in $\Z[I_\alpha]$.  Using the single coordinate $I_\alpha$ on $\sfXA$, we have $\sfXA=\C$.  In particular, we have the following well-known proposition.  

\begin{prop} \label{prop:abelian}
Suppose that $N$ is a knot manifold and $H_1(N;\Z) \cong \Z$.  Let $\sfXA \subseteq \sfX(N)$ be the curve consisting of all abelian representations of $\pie N$.  Then 
\begin{enumerate}
\item $\sfXA$ has exactly one ideal point.  This ideal point detects the unique boundary slope which is trivial in $H_1(N;\Z)$.
\item $T(\sfXA)=\Z[I_\alpha]$.
\end{enumerate}
\end{prop} So, by the proposition, $T(\sfXA)$ is a free $\Z[I_\alpha]$-module with rank one.  

%%%%%%%%%%%%%%%%%%%%%%%%%%%%%%%%%%%

\subsection{Norm curves}  Recall the following definition from \cite{BZ2}.

\begin{dfn}
An algebraic curve component $\sfX$ of $\sfX(N)$ is called a {\it norm curve component} if $I_\gamma \co \sfX \to \C$ is non-constant for every non-trivial peripheral element of $\Gamma$.
\end{dfn}

For example, if $N$ is hyperbolic then there is a hyperbolic curve $\sfX_0 \subseteq \sfX(N)$.  By Proposition 3.1.1 of \cite{CS1}, $\sfX_0$ is a norm curve component.

\begin{prop}\label{prop:norm}
Suppose that $N$ is a knot manifold and $\sfX$ is an irreducible algebraic subset of $\sfX(N)$.  Let $\alpha$ be a slope.  If $\sfX$ is a norm curve component then $\rk_\sfX^\C(\alpha) \geq 2$.  
\end{prop}

\proof Suppose that $\sfX$ is a norm component and that $\rk_X^\C(\alpha)<2$.  By Theorem \ref{thm:main}, $\dim(\sfX)=1$.  Hence, $\rk_\sfX^\C(\alpha)=1$ and we have $f \in \C[X]$ with $f \cdot \C[I_\alpha]=\C[\sfX]$.  In particular, $f$ has an inverse $f^{-1} \in \C[I_\alpha]$.  We have \[ \C[I_\alpha] \ = \ f^{-1} f \cdot \C[I_\alpha] \ =\ f^{-1} \cdot \C[\sfX] \ = \ \C[\sfX].\]  Therefore, $\sfX$ has exactly one ideal point.  On the other hand, \cite[Proposition 4.5]{BZ2} and Theorem \ref{thm:cullershalen} show that $\sfX$ strongly detects at least two boundary slopes, so $\sfX$ must have more than one ideal point. 
\endproof

%%%%%%%%%%%%%%%%%%%%%%%%%%%%%%%%%%%

\subsection{Two-bridge knots}

Suppose that $K$ is a two bridge knot.  Define $N=S^3-K$ and $\Gamma=\pie N$.  Let $\langle \mu, \beta \, | \, \omega \beta = \mu \omega \rangle$ be the standard presentation for $\Gamma$ as given in Section 4.5 of \cite{MR}.   The element $\mu \in \Gamma$ is a meridian for the knot.
 
 Suppose that $\sfX$ is an algebraic component of $\sfX(N)$ which is defined over $\Q$.   Let $x=I_{\mu}$ and $y=I_{\mu \beta}$ in $T(\sfX)$.  It is shown in \cite{HT} that $N$ is small and the slope represented by $\mu$ is not a boundary slope.  By Corollary \ref{cor:irred}, we have a set of polynomials $\{ p_j \}_{j=0}^n$ in  $\Z[x]$ where $\sfX$  is given as the zero set in $\C^2$ of \[ P(x,y) \ = \ y^n - \sum_{j=0}^n p_j(x) y^j\]  and $P$ is irreducible in $\C[x,y]$.  Let $\calB=\{ y^j\}_{j=0}^{n-1}\subset T(\sfX)$.  Since $P$ is zero in $T(\sfX)$, $\calB$  generates $T(\sfX)$ as a $\Z[x]$-module, $T_\Q(\sfX)$ as a $\Q[x]$-module, and $\C[\sfX]$ as a $\C[x]$-module.  In particular, \begin{equation} \label{eq:2bridge}  \rk_\sfX^\Z(\mu) \ \leq \ n.\end{equation}

\begin{prop} \label{prop:2bridge}
$\calB$ is a free basis for $\C[\sfX]$ as a $\C[x]$-module, $T_\Q(\sfX)$ as a $\Q[x]$-module, and $T(\sfX)$ as a $\Z[x]$-module.
\end{prop}

\proof
First, we induct on $k$ to show that $M_k=\langle 1, y, \ldots, y^{k} \rangle_{\C[x]}$ is free of rank $k+1$ whenever $k\leq n-1$.  The base case is immediate since $M_0=\C[x]$.

Since $M_k = M_{k-1} + \langle y^k \rangle$, to see that $M_k = M_{k-1} \oplus \langle y^k \rangle$, it is enough to show that $M_{k-1} \cap \langle y^k \rangle = \{0\}$.  Let $f$ be an element of $M_{k-1} \cap \langle y^k \rangle$.  Then we have $\{ p_j \}_{j=0}^k \subset \C[x]$ such that \[f \ = \ \sum_{j=0}^{k-1} p_j(x) y^j \ = \ p_k(x) y^k \] 
Hence, \[ -p_k(x) y^k + \sum_{j=0}^{k-1} p_j(x) y^j \  \in \ (P) \subset \Q[x,y].\]  But $k < n$ and $P$ is irreducible, so each $p_j$ must be zero.

Therefore, $\calB$ is a free basis for $\C[\sfX]$ as a $\C[x]$-module and using the inequalities (\ref{eq:inequality}) and (\ref{eq:2bridge}) we have \[ \rk_\sfX^\C(\mu) \ = \ \rk_\sfX^\Q(\mu) \ = \rk_\sfX^\Z(\mu) \ = \ n.\]  $\Q[x]$ is a PID, so $\calB$ is a free basis for $T_\Q(\sfX)$ as a $\Q[x]$-module.

To see that $\calB$ is a free basis for $T(\sfX)$ as a $\Z[x]$-module, it remains only to show that this module is free.  To do this, we show that $\langle 1, y, \ldots, y^{k-1} \rangle_{\Z[x]} \cap \langle y^k \rangle_{\Z[x]}=\{0\}$, whenever $k \leq n-1$.  If $f$ is in this intersection, then it is in $M_{k-1} \cap \langle y^k \rangle_{\C[x]}$ which we have already seen is trivial.
\endproof

%%%%%%%%%%%%%%%%%%%%%%%%%%%%%%%%%%%

\begin{example} Let $K$ be a trefoil knot.  A calculation shows that $\sfX(N)=\sfX_A \cup \sfX_0$, where $\sfX_0$ is an irreducible curve.  The polynomial $P$ which defines $\sfX_0$ is $y-1$.  Therefore \[ \rk_{\sfX_0}^\C(\mu) \ = \ \rk_{\sfX_0}^\Q(\mu) \ =\  \rk_{\sfX_0}^\Z(\mu)\ =\ 1.\]  Proposition \ref{prop:norm} reaffirms the well-known fact that $\sfX_0$ is not a norm curve.  
 \end{example}

%%%%%%%%%%%%%%%%%%%%%%%%%%%%%%%%%%%

\begin{example} \label{ex:fig8} Let $K$ be the figure-eight knot.  A calculation shows that $\sfX(N)=\sfX_A \cup \sfX_0$, where $\sfX_0$ is a hyperbolic curve.  The polynomial $P$ which defines $\sfX_0$ is $y^2+(-1-x^2)y+(-1+2x^2)$.  Therefore \[ \rk_{\sfX_0}^\C(\mu) \ = \ \rk_{\sfX_0}^\Q(\mu) \ =\  \rk_{\sfX_0}^\Z(\mu)\ =\ 2.\]  

$H^1(N;\Z_2) \cong \Z_2$ with generator $\sigma$.  $H^1(N;\Z_2)$ acts on $\sfX_0$ by $\sigma(x,y)=(-x,y)$.  It follows that $T(\sfY_0)$ can be identified with the subring (with 1) of $T(\sfX_0)$ generated by $x^2$ and $y$. Moreover, $\sfY_0$ can be identified with the zeros of $y^2+(-1-\xi)y+(-1+2\xi)$ where $\xi=x^2$.  We have \[ \rk_{\sfY_0}^\C(\mu) \ = \ \rk_{\sfY_0}^\Q(\mu) \ =\  \rk_{\sfY_0}^\Z(\mu)\ =\ 2.\] 
 
Let $\lambda=[\beta,\mu^{-1}] [\mu,\beta^{-1}]$.  Then $\lambda$ commutes with $\mu$ and represents the boundary slope determined by a Seifert surface.  Since $\sigma(\lambda)$ is the identity, $l=I_\lambda \in T(\sfY_0)$.  If we take $(l,\xi,y)$ as coordinates on $\sfY_0$ then the ideal determined by $\sfY_0$ is generated by the following polynomials in $\C[l,\xi,y]$.
\begin{align*}
p_1 &=  (-1 - 4 l)+( - 11  + 4 l) y + (16  - l) y^2 - 7 y^3 + y^4\\
p_2 &=  2 - l - 5 \xi + \xi^2\\
p_3 &=  1 - 2 \xi + y + \xi y - y^2\\
p_4 &=  (2l)+\xi+(5-l)y+5y^2+y^3.
\end{align*}
Let $M=\langle 1,y,y^2,y^3 \rangle_{\C[l]}$.  The first polynomial describes the image of $\sfY_0$ under the map $\chi \mapsto (L(\chi),y(\chi))$.  Hence, $p_1$ is irreducible and $M$ is free of rank 4.  The last polynomial shows that $\xi^i y^j \in M$ for every $i, j \in \Z^+$.  It follows that \[ \rk_{\sfY_0}^\C(\lambda) \ = \ \rk_{\sfY_0}^\Q(\lambda) \ =\  \rk_{\sfY_0}^\Z(\lambda)\ =\  4.\]   

Let $s=I_{\mu^2\lambda}$.  Then $s \in T(\sfY_0)$.   Similar calculations show that $\{ 1,y,y^2,y^3,\xi \}$ generates $\C[\sfY_0]$ as a $\C[s]$-module and that \[\rk_{\sfY_0}^\C(\mu^2\lambda) \ = \ \rk_{\sfY_0}^\Q(\mu^2\lambda) \ =\  \rk_{\sfY_0}^\Z(\mu^2\lambda)\ =\  5.\]
\end{example}

%%%%%%%%%%%%%%%%%%%%%%%%%%%%%%%%%%%

\subsection{Other examples}

%%%%%%%%%%%%%%%%%%%%%%%%%%%%%%%%%%%

\begin{example} \label{eg:sister} Let $N$ be the {\tt Snappea} census manifold $M_{003}$ from  \cite{CHW}.  $N$ is a finite volume hyperbolic 3-manifold with a single cusp.  $\langle  \gamma , \eta  \, | \,  \gamma  \eta  \gamma ^{-2} \eta  \gamma  \eta ^3 \rangle$ is a presentation for $\Gamma = \pie N$.  The elements $\mu=( \eta ^2 \gamma  \eta  \gamma )^{-1}$ and $\lambda=( \gamma  \eta  \gamma )^{-1} \eta  \gamma  \eta $ together generate a peripheral subgroup.  Define 
\begin{align*} m=&I_\mu & x&=I_\gamma & y&=I_\eta &  z=&I_{\gamma \eta}
\end{align*}  Let $\sfX_0 \subset \sfX(N)$ be a hyperbolic curve and let $\sfY_0$ be the image of $\sfX_0$ in the $\psl \C$-character variety.  We can take $(m,x,y,z)$ as coordinates on $\sfX_0$.  The ideal in $\C[m,x,y,z]$ determined by $\sfX_0$ is generated by the following three polynomials \begin{align*} p_1&=z^4-mz^2-z^2+1\\ p_2&=-z^2+m+y \\ p_3&=-z^3+mz+z+x.\end{align*}  The first polynomial is irreducible and gives that $\langle 1,z,z^2,z^3 \rangle_{\C[m]}$ is free of rank $4$.  The last two polynomials show that $x, y \in \langle 1,z,z^2,z^3 \rangle_{\C[m]}$.  It follows that $\C[\sfX_0]=\langle 1,z,z^2,z^3 \rangle_{\C[m]}$ and \[ \rk_{\sfX_0}^\C(\mu) \ = \ \rk_{\sfX_0}^\Q(\mu) \ =\  \rk_{\sfX_0}^\Z(\mu)\ =\ 4.\]  

Let $\zeta=z^2$.  It is straightforward to check that $m \in T(\sfY_0)$ and that $(m,\zeta)$ can be taken as coordinates on $\sfY_0$.  As such, $\sfY_0$ is the zero locus of the polynomial $\zeta^2+(-m-1)\zeta+1$.  It follows that 
\[\rk_{\sfY_0}^\C(\mu) \ = \ \rk_{\sfY_0}^\Q(\mu) \ =\  \rk_{\sfY_0}^\Z(\mu)\ =\ 2.\]  We see that these ranks achieve the smallest possible value for the rank of a norm curve.  In contrast to Example \ref{ex:fig8}, $\rk_{\sfY_0}^\C(\mu)<\rk_{\sfX_0}^\C(\mu)$.
\end{example}

%%%%%%%%%%%%%%%%%%%%%%%%%%%%%%%%%%%

\begin{example} Let $N$ be the exterior of the knot $8_{20}$. The knot $8_{20}$ is hyperbolic.  Let $\sfX_0 \subset \sfX(N)$ be a hyperbolic curve.  As outlined in the introduction to this paper, $\Gamma$ is generated by a pair of elements $\mu$ and $\gamma$ where $\mu$ is a meridian.  Define $x=I_\mu$, $y=I_\gamma$, and $z=I_{\gamma \mu}$.  The functions $x,y,z$ give an embedding of $\sfX_0$ into $\C^3$.  

Let $\calG = \{ g_j \}_1^5 \subset \Z[x,y,z]$ where
\begin{align*}
 g_1& = \ 1 + 5 y + 7 y^2 + 2 y^3 - 2 y^4 - y^5 - 2 x^2 - 6 y x^2 - 3 y^2 x^2 + y^4 x^2 + x^4 + y x^4\\
g_2& =\  -x - 3 y x - y^2 x - y^3 x + x^3 + y x^3 + y^2 z + y^3 z\\
g_3& =\   -1 - 4 y - 3 y^2 + y^3 + y^4 + x^2 + 2 y x^2 - y^3 x^2 + y^2 x z\\
g_4& =\   x + 3 y x + 2 y^2 x - x^3 - y x^3 - z - 2 y z - y^2 z + x^2 z\\
g_5& =\   2 + 6 y - 2 y^3 - 3 x^2 - y x^2 + 2 y^2 x^2 + x z - 3 y x z + z^2.
\end{align*}
$\calG$ is a Groebner basis for the ideal $(\calG)$ with respect to the pure lexicographic order on monomials in $\C[x,y,z]$ determined by the relationship $y<x<z$.  Furthermore, $(\calG)$ is the kernel of the natural map $\C[x,y,z] \to \C[\sfX_0]$.

We use the symbol $\equiv$ to indicate congruence modulo $(\calG)$ in $\C[x,y,z]$ and we write $LM(f)$ for the leading monomial of $f \in \C[x,y,z]$ with respect to our chosen monomial order.  For $1 \leq j \leq 5$, let $m_j=LM(g_j)$.  Then $(m_j)_1^5=(yx^4, y^3z, y^2xz, x^2z, z^2)$ and $m_j<m_{j+1}$. 

Oertel's work in \cite{Oe} shows that $N$ is a small knot manifold, so Corollary \ref{cor:irred} implies that there are integral dependencies for $y$ and $z$ over $\Z[x]$ given by irreducible polynomials.  The polynomial $g_1$ is the polynomial for $y$ and  and a quick calculation produces the polynomial  \[ z^5  - 2 x z^4 + (-2 + 3 x^2) z^3 + (12 x - 9 x^3 +x^5) z^2+ (-18 x^2 + 10 x^4 - x^6) z + (6 x^3 - 2 x^5 ) \, \in \, (\calG)\] for $z$.  These polynomials also show that the set $\calB  =  \left\{ y^iz^j \, | \, 0 \leq i,j \leq 4 \right\}$ is a generating set for $\C[\sfX_0]$ as a $\C[x]$-module, $T_\Q(\sfX_0)$ as a $\Q[x]$-module, and $T(\sfX_0)$ as a $\Z[x]$-module.  However, $\calB$ is to large to be a basis for any of these modules.  

\begin{thm} \label{thm:8_20}
The set $\calB'= \{ 1, y, y^2, y^3, z, y^2z \}$ is a free basis for $\C[\sfX_0]$ as a $\C[x]$-module, $T_\Q(\sfX_0)$ as a $\Q[x]$-module, and $T(\sfX_0)$ as a $\Z[x]$-module.  In particular, \[ \rk_{\sfX_0}^\C(\mu) \ = \ \rk_{\sfX_0}^\Q(\mu) \ = \ \rk_{\sfX_0}^\Z(\mu) \ = \ 6.\]
\end{thm}

The proof of Theorem \ref{thm:8_20} follows from the following lemmas.  

\begin{lemma} \label{lem:820a}
If $0 \leq k \leq 4$ then the rank of $\langle 1, y, \ldots, y^k \rangle_{\C[x]}$ is $k+1$.
\end{lemma}
\proof The proof uses the same argument as Proposition \ref{prop:2bridge}. \endproof

\begin{lemma} \label{lem:820b}
If $k \geq 1$ then there is a polynomial $P \in \C[x,y]$ such that \[ x^{2k} z \ \equiv \ P + z(1+2y+y^2)^k.\]
\end{lemma}
\proof Induct on $k$ and use the polynomial $g_4$. \endproof

\begin{lemma} \label{lem:820c}
If $k \geq 3$ then there is a polynomial $P \in \C[x,y]$ such that \[ y^kz \ \equiv \ P \pm y^2z.\]
\end{lemma}
\proof Induct on $k$ and use the polynomial $g_2$. \endproof

\begin{lemma} \label{lem:820d}
If $k \geq 1$ then there is a polynomial $P \in \C[x,y]$ and $d_0, \ldots, d_4 \in \C$ such that \[ x^kz \ \equiv \ P+d_0 z+d_1 yz+d_2 y^2z +d_3 xz +d_4 xyz.\]
\end{lemma}
\proof The lemma is trivial if $k=1$.  

If $k$ is even, a straightforward application of Lemmas \ref{lem:820b} and \ref{lem:820c} gives numbers $d_0, d_1$, and $d_2$ such that $x^kz$ is equivalent to a polynomial in \[d_0z+d_1yz+d_2y^2z+\C[x,y].\]  If $k\geq 3$ is odd, we apply the result in the even case to obtain $c_0, c_1, c_2 \in \C$ where $x^kz$ is equivalent to a polynomial in \[ c_0xz+c_1xyz+c_2xy^2z + \C[x,y].\]  To finish the proof, notice that $g_3$ shows that $xy^2z$ is represented by a polynomial in $\C[x,y]$.  \endproof

\begin{lemma} \label{lem:820e}
Suppose $p \in \C[x,y]$, $g \in \C[x]$, and $d_0, \ldots, d_4 \in \C$. Define \[ f \ = \ p+z\left(g+d_0+d_1y+d_2y^2+d_3x+d_4xy\right).\]  If $f \in (\calG)$ then $f=p$.
\end{lemma}
\proof The proof is by induction on $\deg(g)$.

If $\deg(g)=0$ then $g \in \C$ and \[ f \ = \ p+(g+d_0)z+d_1yz+d_2y^2z+d_3xz+d_4xyz \ \in \ (\calG). \]   Since $\calG$ is a Groebner basis for $(\calG)$, some $m_j$ must divide the leading term of the polynomial on the right.  This forces $f=p$.

Let $k\geq 1$ and assume that the result holds for polynomials with degree less than $k$.  Suppose that $g=\sum_0^k c_j x^j$.   By Lemma \ref{lem:820d}, we have $q \in \C[x,y]$ and $r_0, \ldots, r_4 \in \C$ such that \[ c_kx^kz \ \equiv \ q+r_0z+r_1yz+r_2y^2z+r_3xz+r_4xyz.\]  Therefore, 
\begin{multline*}p+ q + \left( \sum_{j=0}^{k-1} c_j x^j \right)z+(d_0+r_0)z+(d_1+r_1)yz+\\+(d_2+r_2)y^2z+(d_3+r_3)xz+(d_4+r_4)xyz \ \in \ (\calG).
\end{multline*}  By the inductive hypothesis, we have 
\begin{gather*} c_2 = \cdots = c_{k-1} =0\\ c_0+d_0+r_0=0\\c_1+d_3+r_3=0\\ d_1+r_1=d_2+r_2=d_4+r_4=0.
\end{gather*}  Hence, \[ f \ = \ p-r_0z+d_1yz+d_2y^2z-r_3xz+d_4xyz.\]  As before, we must have $f=p$ since $\calG$ is a Groebner basis for $(\calG)$. \endproof

\begin{lemma} \label{lem:820f}
For every $k \in \Z^+$, \[ \langle 1, y, \ldots , y^k, z \rangle_{\C[x]} \ = \ \langle 1, \ldots , y^k \rangle_{\C[x]} \oplus \langle z \rangle_{\C[x]}.\]
\end{lemma}
\proof Take $h \in \langle 1, y, \ldots, y^k \rangle_{\C[x]} \cap \langle z \rangle_{\C[x]}$.  Then there are polynomials $p \in \C[x,y]$ and $g \in \C[x]$ such that $p$ represents $h$ in $\C[\sfX_0]$ and $ p \equiv -gz$.  Define $f=p+gz$ and notice that $f \in (\calG)$.  By Lemma \ref{lem:820e}, $ p \in (\calG)$ and $h=0$. \endproof

\begin{lemma} \label{lem820g}
$\langle \calB' \rangle_{\C[x]}$ is a free $\C[x]$-module with rank 6.
\end{lemma}
\proof Lemmas \ref{lem:820a} and \ref{lem:820f}, reduce the problem to showing that the intersection \begin{equation} \label{intersection} \langle 1, \ldots, y^3, z\rangle_{\C[x]} \, \cap \, \langle y^2z \rangle_{\C[x]}\end{equation} is trivial.

Since $g_3 \in \calG$, we have that $\langle 1, \ldots, y^3, xy^2z \rangle_{\C[x]} = \langle 1, \ldots, y^4 \rangle_{\C[x]}$.  Together with Lemmas \ref{lem:820a} and \ref{lem:820f}, this shows that $\langle 1, \ldots, y^3, xy^2z, z \rangle_{\C[x]}$ is free with rank six.  

Now suppose that $h$ is an element of the intersection (\ref{intersection}).  Then $xh$ is an element of both $\langle 1, \ldots, y^3, z\rangle_{\C[x]}$ and $\langle xy^2z \rangle_{\C[x]}$.  But we know that the intersection of this module is trivial.  Hence, $hx=0$ and so $h=0$.
\endproof

Now to prove Theorem \ref{thm:8_20}, we need only show that $\calB'$ spans the modules in question.  This is straightforward using the polynomials $\calG$.
\end{example}

%%%%%%%%%%%%%%%%%%%%%%%%%%%%%%%%%%%

\begin{example} 
$N$ is the once punctured hyperbolic torus bundle from Section 5 of \cite{D2}.  For more details in the calculations that follow, see \cite{D2}.  We have
\[ \Gamma = \left\langle \alpha,\beta,\tau \, \big| \, \tau\alpha\tau^{-1}=(\beta\alpha\beta)^{-1}, \tau\beta\tau^{-1}=\beta\alpha(\beta\alpha\beta)^{-3} \right\rangle.\]  The elements $\tau$ and $\lambda=[\alpha,\beta]$ form a basis for a peripheral subgroup of $N$.  The functions 
\begin{align*} 
  t&=I_\tau &  u&=I_{\alpha\tau} & v&=I_{\beta\tau}  &  w&=I_{\alpha\beta\tau} \\
  x&=I_\alpha & y&=I_\beta & z&=I_{\alpha\beta}  & l&=I_\lambda
\end{align*}
give an embedding of $\sfX(N)$ into $\C^8$.  $\sfX(N)$ contains exactly two hyperbolic curves, we label them $\sfX_\epsilon$ where $\epsilon \in \{ 0,1 \}$.  It is straightforward to verify that $\sfX_\epsilon$ strongly detects the slopes represented by $t$ and $lt$ and does not detect a closed essential surface.  This is done in \cite{D2}.  It is also easy to use Theorem \ref{thm:I_mu} to confirm these facts.

\begin{thm} \label{thm:bundle} Let $\epsilon \in \{ 0, 1\}$.
%\[\] \vspace{-.35in}
\begin{enumerate}
\item $\rk_{\sfX_\epsilon}^\C(\lambda) = 8$
\item $\rk_{\sfX_\epsilon}^\Q(\lambda) = 16$, and 
\item $\rk_{\sfX_\epsilon}^\Z(\lambda) = 17$.
\item As a $\Z[l]$-module, $T(\sfX_\epsilon)$ is torsion free but not free.
\end{enumerate}
\end{thm}

We prove Theorem \ref{thm:bundle} as a series of lemmas.  

Suppose $k$ is a field, let $R=k[x_0, \ldots, x_n]$, and $\calI \subseteq R$ and ideal.  Take $<$ to be the monomial order generated by $x_i < x_j$ iff $i<j$ and let $\calG=\{ g_j\}_1^r$ be a Groebner basis for $\calI$ with respect to $<$.  Write $m_j = LT(g_j)$.  For $f \in R$, write $\bar{f}$ to denote its image in $R/\calI$.

\begin{lemma} \label{lem:free}
Suppose that $m_j \in k[x_1, \ldots x_n]$ for every $j$.  If $\{ n_j \}_1^k$ is a collection of monomials in $k[x_1, \ldots , x_n]$ such that 
\begin{enumerate}
\item The $k[x_0]$-submodule of $R/\calI$ generated by $\{ \bar{n}_j \}_1^{k-1}$ is free of rank $k-1$,
\item $n_i < n_j$ for every $i<j$, and
\item $n_k$ is not divisible by any $m_j$
\end{enumerate}
then the $k[x_0]$-submodule generated by $\{ \bar{n}_j \}_1^k$ is free of rank $k$.
\end{lemma}

\proof
Suppose $f \in R$ and \[\bar{f} \in \langle \bar{n}_1, \ldots , \bar{n}_{k-1} \rangle_{k[x_0]} \cap \langle \bar{n}_k \rangle_{k[x_0]}.\]  Then there is a set of polynomials $\{ p_j\}_1^k \subset k[x_0]$ such that \[ \bar{f} \ = \ -p_k \bar{n}_k \ =\ \sum_{j=1}^{k-1} p_j \bar{n}_j.\]  In particular \[ g \ := \ \sum_{j=1}^k p_j n_j \ \in \ \calI.\]  If $p_k=0$ then $\bar{f}=0$ and we are done.  Otherwise, let $cx_0^s$ be its leading term.  Then, by (2), $LT(g) = cx_0^s n_k$.  Since $g \in \calI$ and $\calG$ is a Groebner basis for $\calI$, we must have that $LT(g)$ is divisible by some $m_j$.  This contradicts (3), so we must have $p_k=0$.
\endproof

Let $\calI$ be the kernel of the natural map $\C[l,t,u,v,w,x,y,z] \to \C[\sfX_\epsilon]$.  The following polynomials are in $\calI$
\begin{align*}
& 8u+tyz-ty^3z & & 4v-6ty+ty^3 \\ &4w+2tz-ty^2z && 2x-yz.
\end{align*}
\begin{remark}
Recall that Groebner bases can be used to solve the ideal membership problem. \end{remark}

This shows that the projection given by $(l,t,u,v,w,x,y,z) \mapsto (l,t,y,z)$, restricts to an isomorphism (defined over $\Q$) from $\sfX_0 \cup \sfX_1$ onto its image.  Hence, we 
have natural identifications of $\C[\sfX_\epsilon]$ and $\C[\sfX_0 \cup \sfX_1]$ with quotients of $\C[l,t,y,z]$.  Define 
\begin{align*}
g_1 & = 16+(-4-2l)z^2+z^4 \\
g_2 & = (-8-2l)+z^2+2y^2\\
g_3 &= 2t+(-1)^{\epsilon}  \cdot i z^2 \\
h &= -8+(2+l)z^2+2t^2
\end{align*}
and $\calG=\{ g_1, g_2, g_3\}$.  Let $<$ be the pure lexicographic monomial order for $\C[l,t,y,z]$ induced by the relationship $l<z<y<t$.  It can be easily verified that $\calG$ is a Groebner basis for $(\calG)$ and $(\calG)$ is the kernel of the natural map $\C[l,t,y,z] \to \C[\sfX_\epsilon]$.  Similarly, $\{ g_1, g_2, h\}$ is a Groebner basis for the ideal it generates and this ideal is the kernel of the natural map $\C[l,t,y,z] \to \C[\sfX_0 \cup \sfX_1]$. 

\begin{lemma} \label{lem:bundleC}
 The set \[\calB_\C\ =\ \{ 1,z,z^2,z^3,y,zy, z^2y, z^3y\}\] is a basis for $\C[\sfX_\epsilon]$ as a $\C[l]$-module.  \end{lemma}
\proof
Since the natural map $\C[l,t,y,z]/(\calG) \to \C[\sfX_\epsilon]$ is an isomorphism, we need only show that $\calB_\C$ is a $\C[l]$-module basis for $\C[l,t,y,z]/(\calG)$.  

Let $\calM \subseteq \C[l,t,y,z]/(\calG)$ be the $\C[l]$-submodule generated by $\calB_\C$.  Lemma \ref{lem:free} shows that $\calM$ is free with rank $8$.  

The polynomial $g_3$ shows that $t$ may be expressed in terms of $z$.  So, to finish the proof, it is enough to show that $y^j z^k \in \calM$ for every pair of integers $j,k \geq 0$.  By definition of $\calB_\C$ and $g_1$, $z^j \in \calM$, for every integer $j \geq 0$.  Using $g_2$, we also have that $y^2z^j \in \calM$, for every integer $j \geq 0$.  Now, using this fact and the definition of $\calB_\C$, we conclude that if $f \in \calM$ then $yf \in \calM$.  In particular, $y^j z^k \in \calM$ for every pair of integers $j,k \geq 0$.  Therefore, $\calM = \C[l,y,z]/(\calG)$.
\endproof

\begin{lemma} \label{lem:bundleQ}
The set \[ \calB_\Q\ =\  \calB_\C \cup t \calB_\C\] is a basis for $T_\Q(\sfX_\epsilon)$ as a $\Q[l]$-module. 
\end{lemma}
\proof  The natural map $\C[l,t,y,z] \to \C[\sfX_\epsilon]$ restricts to a surjection $\Q[l,t,y,z] \to T_\Q(\sfX_\epsilon)$.  The kernel of this map is $(\calG) \cap \Q[l,t,y,z]$.  We claim that $\{ g_1, g_2, h\}$ is a Groebner basis for this ideal with respect to the monomial order $<$.  It is immediate to verify that this is a Groebner basis for $(g_1, g_2, h)$, so we need only show that $(g_1, g_2, h) = (g_1, g_2, g_3) \cap \Q[l,t,y,z]$.  The inclusion $(g_1, g_2, h) \subseteq (g_1, g_2, g_3) \cap \Q[l,t,y,z]$ is clear.  To prove the opposite, suppose $f \in (\calG) \cap \Q[l,t,y,z]$.  Then  $f(\chi)=0$ for every $\chi \in \sfX_\epsilon$.  Observe that complex conjugation induces a bijection $\sfX_0 \to \sfX_1$.  So $f(\bar{\chi})$ is also zero.  In particular, $f$ is zero on $\sfX_0 \cup \sfX_1$.  Therefore $f \in (g_1, g_2, h)$.

Let $\calM \subset T_\Q(\sfX_\epsilon)$ be the $\Q[l]$-submodule generated by $\calB_\Q$.  Since $\{ g_1, g_2, h \}$ is a Groebner basis for the kernel of the map $\Q[l,t,y,z] \to T_\Q(\sfX_\epsilon)$, Lemma \ref{lem:free} shows that $\calM$ is free with rank $16$.  The remainder of the proof follows an argument similar to the last part of the proof of Lemma \ref{lem:bundleC}.
\endproof

Let \[ \calB_\Z\ =\ \{ 1, y, y^2, y^3, t, u, v, w, x, z, ty, uy, xy, xy^2, tx, ux, vy \}. \]

\begin{lemma} \label{lem:Qbasis}
$\calB_\Z - \{ vy \}$ is a basis for $T_\Q(\sfX_\epsilon)$ as a $\Q[l]$-module.
\end{lemma}
\proof
The following equations hold in $\C[\sfX_\epsilon]$.
\begin{align*}
z^2&=2((4+l)-y^2) & z^3&=2((4+l)z-2xy) \\
yz&=2x & yz^2&=2((4+l)y-y^3) \\
yz^3&=2((8+2l)x-2xy^2) & tz&=(l+4)w-2uy\\
tz^2&=2((l+4)t-2ux) & tz^3&=2((l^2+6l+4)w-2(l+2)uy)\\
tyz&=2tx  & tyz^2&=2((l-2)ty+4v)  \\
 tyz^3&=4((l+2)tx-4u)
\end{align*}
So, every element of $\calB_\Q$ can be expressed as a $\Q[l]$-linear combination of the elements from $\calB_\Z-\{vy\}$.  Since the rank of this module equals $| \calB_\Z-\{ vy\}|$, this set is a basis.
\endproof

\begin{lemma} \label{lem:Zgenerates}
$\calB_\Z$ generates $T(\sfX_\epsilon)$ as a $\Z[l]$-module.
\end{lemma}
\proof
The following equations hold in $\C[\sfX_\epsilon]$.  It is tedious but straightforward to use them to show that $\calB_\Z$ generates $T(\sfX_\epsilon)$ as a $\Z[l]$-module.
\begin{align*}
y^2z&=2yx   & x^2u&=2u+xt & uy^2x&=(6+l)xu-(6+l)t  \\
 yz&=2x & x^2&=-y^2+(6+l)  & yxt&=(6+l)w-2uy\\ 
yw&=2u &  y^4&=(6+l)(y^2-2) & zy^2x&=2(-y^3+(6+l)y) \\
 zyu&=2xu & zxu&=(6+l)w-2uy & zyx&=2(-y^2+(6+l)) \\
y^3z&=2y^2x  & zx&=(4+l)y-y^3 & vyu&=3y^2x-2(6+l)x \\ 
zu&=yt & u^2&=(-6-l)+y^2 & vy^2x&=(6+l)(3w-yu)    \\
xw&=yt & tyu&=2xy-(6+l)z & wy^3&=2((4+l)u-xt) \\ 
\end{align*}

\vspace{-.2in}
\begin{align*}
t^2&=(-4-6l-l^2)+(2+l)y^2 & t^2x&=(-4-6l-l^2)x+(2+l)y^2x \\                                                               
t^2y&=(-4-6l-l^2)y+(2+l)y^3 & txu&=(-12-8l-l^2)+(4+l)y^2\\
uxt&=(-12-8l-l^2)+(4+l)y^2 & vxt&=2(l+3)yx-(l+6)(l+1)z \\
wxt&=(-4-6l-l^2)y+(2+l)y^3 & vyt&=2((l+3)y^2-(l+1)(l+6))  \\
vxu&=(-6-7l-l^2)y+(3+l)y^3 &    tv&=(l+1)y^3-(5l+l^2)y\\
v^2&=(-12l-2l^2)+(3+2l)y^2 
\end{align*}
\begin{align*}
zw&=2t & y^2u&=(4+l)u-xt & x^2y&=-y^3+(6+l)y\\ 
wy^2&=2uy & y^3z&=(6+l)(yx-z)  & x^2t&=-2xu+(6+l)t \\ 
wyx&=2ux & xv&=(6+l)w-3uy &  uy^3&=(6+l)(yu-w)\\
xyu&=3ty-2v  & uv&=3yx-(6+l)z & uw&=-(4+l)y+y^3 \\
x^2y^2&=2(6+l)   & tw&=2yx-(4+l)z & wyu&=2(y^2-(6+l)) \\
ty^2x&=4u+2xt & wxu&=2yx-(6+l)z\ &  wyt&=2(y^2x-(4+l)x)\\ 
wy^2x&=6yt-4v & u^2y&=(-6-l)y+y^3 & u^2x&=(-6-l)x+y^2x\\
vyx&=-lu+3xt & 2vy&=(6+l)t-lxu & tu&=y^2x-(4+l)x\\
yxu&=3ty-2v & vy^2&=(3-l)yt+lv & w^2&=2((-4-l)+y^2)  \\
zv&=2xt-2u & zxt&=(l-2)yt+4v  & wv&=2y^2x-2(3+l)x \\
\end{align*} 

\vspace{-.45in}
\endproof

\begin{lemma} \label{lem:directsum}
Define \[ \calM_1 \ =\ \Big\langle \calB_\Z-\{t,ux,vy\} \Big\rangle_{\Z[l]}\subset T(\sfX_\epsilon) \quad \text{and} \quad \calM_2 \ =\ \Big\langle t, ux, vy \Big\rangle_{\Z[l]} \subset T(\sfX_\epsilon).\]
Then $ T(\sfX_\epsilon)  = \calM_1 \oplus \calM_2$.
\end{lemma}
\proof Lemma \ref{lem:Zgenerates} shows that $T(\sfX_\epsilon) = \calM_1 + \calM_2$, so we need only show that $\calM_1 \cap \calM_2 = \{ 0 \}$.  

Suppose $f \in \calM_1 \cap \calM_2$ and label the elements of $\calB_\Z-\{t,ux,vy\}$ as $\{ b_j \}_1^{14}$.  Then we have $\{ p_j \}_1^{14} \cup \{ q_j \}_1^3 \subset \Z[l]$ with
\begin{equation*} %\label{eqn:directsum}
f \ = \ \sum_{j=1}^{14} p_j b_j \ = \ q_1 t + q_2 ux + q_3 vy.
\end{equation*}  
This yields the expression \[ 2 q_3 vy \ =\ -2q_1 t-2q_2 ux+2 \sum_{j=1}^{14} p_j b_j.\]
However, $\calB_\Z-\{vy\}$ is a $\Q[l]$-basis for $T_\Q(\sfX_\epsilon)$, so the only way to write $2 vy$ as a $\Z[l]$-linear combination in $\calB_\Z-\{vy\}$ is listed in the proof of Lemma \ref{lem:Zgenerates}, namely \[ 2 vy \ =\ (6+l)t-lux.\]  Hence $p_j=0$ for every $1\leq j \leq 14$.
\endproof

\begin{lemma} \label{lem:projective}
$\calM_2$ is not projective. 
\end{lemma}
\proof Let \[\pi \co \bigoplus_1^3 \Z[l] \to \calM_2 \] be the epimorphism given by $(p,q,r) \mapsto pt+qux+rvy$.  Argue as in the proof of Lemma \ref{lem:directsum} to see that $\ker(\pi)$ is generated by $(l+6, -l, -2)$.  Let $\calJ$ be the ideal $(2,l) \subset \Z[l]$ and view $\calJ$ as a $\Z[l]$-module.  The map \[ \eta \co \bigoplus_1^3 \Z[l] \to \calJ\] given by $(p,q,r) \mapsto -2q+lr$ is also an epimorphism.  We have $\ker(\pi) \subseteq \ker(\eta)$, so there is an epimorphism \[ \gamma \co \calM_2 \to \calJ\] with $\gamma \pi = \eta$.  In particular, 
\begin{align*}
\gamma(t)&=0 & \gamma(ux)&=-2 & \gamma(vy)&=l.
\end{align*}
Let \[ \phi \co \Z[l] \oplus \Z[l] \to \calJ \] be the epimorphism given by $(p,q) \mapsto 2p+lq$.

Suppose for a contradiction, that $\calM_2$ is projective.  Then $\gamma$ lifts to an epimorphism \[ \psi \co \calM_2 \to \Z[l] \oplus \Z[l].\]  For $1 \leq j \leq 3$, we have $(p_j, q_j) \in \Z[l] \oplus \Z[l]$ with \begin{align*}
\psi(t)&=(p_1,q_1) & \psi(ux)&=(p_2,q_2) & \psi(vy)&=(p_3,q_3).
\end{align*}
Then
\begin{align*}
0&=\gamma(t)=\phi(p_1,q_1)=2p_1+lq_1\\
-2&=\gamma(ux)=\phi(p_2,q_2)=2p_2+lq_2\\
l&=\gamma(vy)=\phi(p_3,q_3)=2p_3+lq_3.
\end{align*}
Notice that the constant terms of $p_1$ and $p_3$ must be zero and the constant term of $p_2$ is $-1$.  Apply $\psi$ to  the equation $2vy+lux=(6+l)t$ to see that the equation \[2p_3+lp_2=(6+l)p_1\] must hold in $\Z[l]$.  For a contradiction, observe that the the degree one coefficient of the left hand side of this equation is odd and on the right hand side the coefficient is even.
\endproof

\begin{lemma} \label{lem:Zbasis}
The $\Z[l]$-module $T(\sfX_\epsilon)$ is torsion free but not free and $\rk_{\sfX_\epsilon}^\Z(\lambda)=17$.
\end{lemma}
\proof $\C[\sfX_\epsilon]$ is an integral domain, so $T(\sfX_\epsilon)$ is torsion free as a $\Z[l]$-module.  Direct summands of free modules are projective. Therefore, by Lemmas \ref{lem:directsum} and \ref{lem:projective}, $T(\sfX_\epsilon)$ is not a free $\Z[l]$-module.

By Lemma \ref{lem:Zgenerates}, $\rk_{\sfX_\epsilon}^\Z(\lambda) \leq 17$.  Suppose that $\{ a_j \}_1^{16}$ generates $T(\sfX_\epsilon)$ as a $\Z[l]$-module.  Then every element of $\calB_\Z$ can be written as a $\Z[l]$-linear combination of the $a_j$'s.  By Lemma \ref{lem:bundleQ}, $\{ a_j\}_1^{16}$ is a free $\Q[l]$-basis for $T_\Q(\sfX_\epsilon)$. It follows that the $\Z[l]$-module generated by $\{ a_j \}_1^{16}$ is also free.  This is a contradiction because we know that $T(\sfX_\epsilon)$ is not a free $\Z[l]$-module.  Therefore, $\rk_{\sfX_\epsilon}^\Z(\lambda)=17$.
\endproof

\end{example}

\bibliographystyle{plain}

\bibliography{holes}

\end{document}